\newif\ifTR
\global\TRtrue

\ifTR
 \documentclass{article}
 \usepackage{natbib}%
 \usepackage{hyperref}  
\usepackage{orcidlink} 
\usepackage{authblk}

\newif\ifbackmatter%
\newcommand{\backmatter}{\global\backmattertrue}%

\usepackage[title]{appendix}%

 \else
 \documentclass[pdflatex,sn-basic]{sn-jnl}
 \fi

\usepackage{epstopdf}
\usepackage{graphicx}%
\usepackage{multirow}%
\usepackage{amsmath,amssymb,amsfonts}%
\usepackage{amsthm}%
\usepackage{mathrsfs}%
\usepackage[title]{appendix}%
\usepackage{xcolor}%
\usepackage{textcomp}%
\usepackage{manyfoot}%
\usepackage{booktabs}%
\usepackage{algorithm}%
\usepackage{algorithmicx}%
\usepackage{algpseudocode}%
\usepackage{listings}%

\usepackage{multirow}%
\usepackage{amsmath,amssymb,amsfonts}%
\usepackage{amsthm}%
\usepackage{mathrsfs}%
\usepackage[title]{appendix}%
\usepackage{xcolor}%
\usepackage{textcomp}%
\usepackage{manyfoot}%
\usepackage{booktabs}%
\usepackage{algorithm}%
\usepackage{algorithmicx}%
\usepackage{algpseudocode}%
\usepackage{listings}%

\usepackage{times}
\usepackage{pgfplotstable}
\usepackage{comment}

\usepackage{amsmath,amsfonts,bm}

\newcommand{\ML}{machine learning}
\newcommand{\NN}{neural network}
\newcommand{\NNs}{neural networks}
\newcommand{\RNN}{Recurrent Neural Network}

\newcommand{\LSTM}{Long Short-Term Memory}

\newcommand{\MLPs}{Multi-Layer Perceptrons}

\newcommand{\LR}{Lagrangian relaxation}
\newcommand{\LP}{Lagrangian subproblem}

\newcommand{\Lp}{Linear problem}
\newcommand{\LDP}{Lagrangian Dual problem}
\newcommand{\LD}{Lagrangian Dual}

\newcommand{\MP}{Master problem}
\newcommand{\DMP}{Dual Master problem}

\newcommand{\LM}{Lagrangian multiplier}
\newcommand{\LMs}{Lagrangian multipliers}
\newcommand{\CR}{Continuous relaxation}

\newcommand{\BM}{Bundle method}
\newcommand{\BMs}{Bundle methods}
\newcommand{\ABM}{Aggregated Bundle method}

\newcommand{\PBMs}{Proximal Bundle methods}

\newcommand{\CP}{Cutting Plane}

\newcommand{\CPM}{Cutting Plane method}
\newcommand{\SGM}{subgradient method}

\def\1{\bm{1}}

\def\vtheta{{\bm{\theta}}}

\def\valpha{{\bm{\alpha}}}

\def\vb{{\bm{b}}}

\def\vd{{\bm{d}}}

\def\vg{{\bm{g}}}
\def\vh{{\bm{h}}}

\def\vk{{\bm{k}}}

\def\vq{{\bm{q}}}
\def\vr{{\bm{r}}}

\def\vw{{\bm{w}}}
\def\vx{{\bm{x}}}

\def\vpi{{\bm{\pi}}}
\def\vgamma{{\bm{\gamma}}}

\def\vtheta{{\bm{\theta}}}

\def\mA{{\bm{A}}}

\def\mC{{\bm{C}}}

\def\mW{{\bm{W}}}

\DeclareMathAlphabet{\mathsfit}{\encodingdefault}{\sfdefault}{m}{sl}
\SetMathAlphabet{\mathsfit}{bold}{\encodingdefault}{\sfdefault}{bx}{n}

\def\sR{{\mathbb{R}}}

\def\vtheta{{\bm{\theta}}}

\def\vb{{\bm{b}}}

\def\vd{{\bm{d}}}

\def\vg{{\bm{g}}}
\def\vh{{\bm{h}}}

\def\vk{{\bm{k}}}

\def\vq{{\bm{q}}}
\def\vr{{\bm{r}}}

\def\vw{{\bm{w}}}
\def\vx{{\bm{x}}}

\def\mA{{\bm{A}}}

\def\mC{{\bm{C}}}

\def\mW{{\bm{W}}}

\def\sR{{\mathbb{R}}}

\usepackage{hyperref}
\usepackage{url}
\usepackage{multirow}
\usepackage{booktabs}
\usepackage{amsmath,amsfonts,amssymb,amsthm}
\usepackage{bm}
\usepackage{eso-pic}
\usepackage{xspace}
\usepackage{amsmath}
\usepackage{svg}
\usepackage{makecell}
\usepackage{microtype}
\usepackage{subfigure}
\usepackage{algpseudocode}
\usepackage{algorithm}
\usepackage{amsmath}
\usepackage{amssymb}
\usepackage{mathtools}
\usepackage{amsthm}
\usepackage[capitalize,noabbrev]{cleveref}
\usepackage{tikz}

\usetikzlibrary{backgrounds,calc,shadows}
\tikzstyle{block} = [rectangle, draw, rounded corners, minimum height=1cm, minimum width=2cm, text centered, thick]
\tikzstyle{embedding} = [block, fill=red!20]
\tikzstyle{positional} = [circle, draw, minimum size=7mm, thick]
\tikzstyle{attention} = [block, fill=orange!20]
\tikzstyle{norm} = [block, fill=yellow!20]
\tikzstyle{feedforward} = [block, fill=blue!20]
\tikzstyle{arrow} = [thick,->,>=Stealth]

\usetikzlibrary{arrows}
\usetikzlibrary{automata}
\usetikzlibrary{shapes.geometric, arrows, chains, fit, positioning}

\tikzstyle{var_node} = [circle, minimum width=1cm, minimum height=1cm, text centered, draw=black, fill=red!30]

\tikzstyle{constr_node} = [circle, minimum width=1cm, minimum height=1cm, text centered, draw=black, fill=blue!30]

\tikzstyle{process} = [rectangle, minimum width=1.5cm, minimum height=1cm,rounded corners, text centered, draw=black,text width=1.8cm, fill=yellow!10]

\tikzstyle{data} = [rectangle, minimum width=2cm, minimum height=1cm, rounded corners,text centered, draw=black, text width=2cm,fill=cyan!05]

\tikzstyle{d_data} = [rectangle,  minimum width=3.2cm, minimum height=0.5cm, rounded corners,text centered, draw=black, text width=2cm, fill=cyan!05]

\tikzstyle{mini_data} = [rectangle, minimum width=.5cm, minimum height=1.5cm, rounded corners,text centered, draw=black, text width=.5cm, fill=white!50]

\tikzstyle{huge} = [rectangle,  minimum width=17.5cm, minimum height=5cm, rounded corners,text centered, draw=black, text width=2cm, fill=cyan!05]

\tikzstyle{large} = [rectangle,  minimum width=7cm, minimum height=3.5cm, rounded corners,text centered, draw=black, text width=2cm, fill=cyan!05]

\tikzstyle{medium} = [rectangle,  minimum width=1cm, minimum height=2.5cm, rounded corners,text centered, draw=black, text width=1.5cm, fill=cyan!05]

\tikzstyle{small} = [rectangle,  minimum width=01.4cm, minimum height=0.5cm,text centered, draw=black, text width=1cm, fill=cyan!05]

\tikzstyle{medium_small} = [rectangle,  minimum width=1.4cm, minimum height=.3cm,text centered, draw=black, text width=1cm, fill=white!05]

\tikzstyle{dot} = [circle, minimum width=0.3cm, minimum height=0.3cm, rounded corners,text centered, draw=black,  fill=white!50]

\tikzstyle{blackbox} = [rectangle,  minimum width=3.5cm, minimum height=5.5cm, rounded corners,text centered, draw=black, text width=2cm, fill=orange!30]

\tikzstyle{arrow} = [thick,->,>=stealth, line width=1.5pt]
\tikzstyle{und_arrow} = [thick,--]

\tikzstyle{c_data} = [circle, minimum width=1cm,text centered, draw=black, text width=1.5cm,fill=cyan!15]
\tikzstyle{c_action} = [rectangle,text width=1.5cm,minimum height=1.5cm,minimum width=1cm, draw=black]
\tikzstyle{big_action} = [rectangle,text width=7cm,minimum height=1.5cm,minimum width=1cm, draw=black]
\tikzstyle{r_flow} = [rectangle,  text width=2.5cm, minimum width=3cm,minimum height=5cm, fill=yellow!15, text centered, draw=black, rounded corners]
\tikzstyle{rr_flow} = [rectangle,  text width=2.5cm, minimum width=3cm,minimum height=5cm, fill=green!15, text centered, draw=black]
\tikzstyle{z_tiny} = [rectangle,  text width=0.25cm, minimum width=0.25cm,minimum height=0.25cm, fill=white!15, text centered, draw=black]
\tikzstyle{r_mlp} = [rectangle,  text width=1.5cm, minimum width=0.75cm,minimum height=0.5cm, fill=white!15, text centered, draw=black]
\tikzstyle{pi_tiny} = [rectangle,  text width=0.5cm, minimum width=0.5cm,minimum height=0.5cm, fill=white!15, text centered, draw=black]

\pgfplotsset{compat=1.18}

\makeatletter
\DeclareRobustCommand\onedot{\futurelet\@let@token\@onedot}
\def\@onedot{\ifx\@let@token.\else.\null\fi\xspace}

\ifTR
\title{Bundle Network: a Machine Learning-Based Bundle Method}

\author[1,2]{Francesca Demelas\orcidlink{0000-0003-1888-3182}\thanks{demelas@lipn.univ-paris13.fr}}
\author[1]{Antonio Frangioni\orcidlink{0000-0002-5704-3170}\thanks{antonio.frangioni@unipi.it}}
\author[2]{Mathieu Lacroix\orcidlink{0000-0001-8385-3890}\thanks{lacroix@lipn.univ-paris13.fr}}
\author[2]{Joseph Le Roux\orcidlink{0000-0002-3889-8536}\thanks{leroux@lipn.univ-paris13.fr}}
\author[3]{Emiliano Traversi\orcidlink{0000-0003-4673-3982}\thanks{emiliano.traversi@essec.edu}}
\author[2]{Roberto Wolfler Calvo\orcidlink{0000-0002-5459-5797}\thanks{wolfler@lipn.univ-paris13.fr}}

\affil[1]{Department of Computer Science, Università di Pisa, 3 Largo B. Pontecorvo, 56127 Pisa, Italy}
\affil[2]{Laboratoire d'Informatique de Paris Nord, Université Sorbonne Paris Nord, 99 avenue Jean-Baptiste Clément, 93430 Villetaneuse, France}
\affil[3]{Department of Information Systems, Data Analytics and Operations, ESSEC Business School, 3 avenue Bernard Hirsch, 95000 Cergy, France}

\date{} 

\else
\title[Bundle Network]{Bundle Network: a Machine Learning-Based Bundle Method}

\author*[1,2]{\fnm{Francesca} \sur{Demelas}\orcid{https://orcid.org/0000-0003-1888-3182}}\email{demelas@lipn.univ-paris13.fr}

\author[1]{\fnm{Antonio} \sur{Frangioni}\orcid{https://orcid.org/0000-0002-5704-3170}}\email{antonio.frangioni@unipi.it}

\author[2]{\fnm{Mathieu} \sur{Lacroix}\orcid{https://orcid.org/0000-0001-8385-3890}}\email{lacroix@lipn.univ-paris13.fr}

\author[2]{\fnm{Joseph} \sur{Le~Roux}\orcid{https://orcid.org/0000-0002-3889-8536}}\email{leroux@lipn.univ-paris13.fr}

\author[3]{\fnm{Emiliano} \sur{Traversi}\orcid{https://orcid.org/0000-0003-4673-3982}}\email{emiliano.traversi@essec.edu}

\author[2]{\fnm{Roberto} \sur{Wolfler~Calvo}\orcid{https://orcid.org/0000-0002-5459-5797}}\email{wolfler@lipn.univ-paris13.fr}

\affil*[1]{\orgdiv{Department of Computer Science}, \orgname{Università di Pisa}, \orgaddress{\street{3 Largo B. Pontecorvo}, \city{Pisa}, \postcode{56127}, \state{Italy}}}

\affil[2]{\orgdiv{Laboratoire d'Informatique de Paris Nord}, \orgname{Université Sorbonne Paris Nord}, \orgaddress{\street{99 avenue Jean-Baptiste Clément}, \city{Villetaneuse}, \postcode{93430}, \state{France}}}

\affil[3]{\orgdiv{Department of Information Systems, Data Analytics and Operations}, \orgname{ESSEC Business School}, \orgaddress{\street{3 avenue Bernand Hirsch}, \city{Cergy}, \postcode{95000}, \state{France}}}

\abstract{
  This paper presents \textit{Bundle Network}, a learning-based algorithm inspired by the Bundle Method for convex non-smooth minimization problems.
  Unlike classical approaches that rely on heuristic tuning of a regularization parameter, our method automatically learns to adjust it from data.
  Furthermore, we replace the iterative resolution of the optimization problem that provides the search direction ---traditionally computed as a convex combination of gradients at visited points--- with a recurrent neural model equipped with an attention mechanism.
  By leveraging the unrolled graph of computation, our \textit{Bundle Network} can be trained end-to-end via automatic differentiation.
  Experiments on Lagrangian dual relaxations of the Multi-Commodity Network Design and Generalized Assignment problems demonstrate that our approach consistently outperforms traditional methods relying on grid search for parameter tuning, while generalizing effectively across datasets.}

\fi

\begin{document}

\ifTR

\maketitle
\abstract{
  This paper presents \textit{Bundle Network}, a learning-based algorithm inspired by the Bundle Method for convex non-smooth minimization problems.
  Unlike classical approaches that rely on heuristic tuning of a regularization parameter, our method automatically learns to adjust it from data.
  Furthermore, we replace the iterative resolution of the optimization problem that provides the search direction ---traditionally computed as a convex combination of gradients at visited points--- with a recurrent neural model equipped with an attention mechanism.
  By leveraging the unrolled graph of computation, our \textit{Bundle Network} can be trained end-to-end via automatic differentiation.
  Experiments on Lagrangian dual relaxations of the Multi-Commodity Network Design and Generalized Assignment problems demonstrate that our approach consistently outperforms traditional methods relying on grid search for parameter tuning, while generalizing effectively across datasets.\\
  \vspace{0.25cm}
  }
\else
\maketitle
\fi

\section{Introduction}

This work presents a \ML{}-based approach to tackle convex non-smooth minimization problems, inspired by bundle-type methods~\citep{lemarechal1981bundle,lemarechal1995new}, i.e., iterative algorithms that generalize the \SGM{} in that the search direction is a convex combination of the gradients at the visited points obtained by solving an optimization problem, called \MP{}. At each iteration, a crucial scalar parameter related to the \emph{regularization} in the \MP{} must be properly chosen to obtain optimal performance. Although there are many variants of the \BM{}~\citep{frangioni2020}, in this work we only consider the best known \emph{proximal} version where the \textit{regularization parameter}, denoted by $\eta$, is the weight of the quadratic \emph{regularization term} in the objective of the \MP{}, weighting the distance of the new point from the better point found so far. This parameter significantly influences the next iterate produced at each iteration in two different ways:
\begin{enumerate}
 \item by changing the \emph{convex multipliers} that are used to form the search direction out of the currently available \emph{bundle of} (approximate) \emph{subgradient information};
 \item by acting as the stepsize along said direction.
\end{enumerate}
Traditionally, selecting this parameter relies on heuristic approaches, referred to as \textit{$\eta$-strategies}, often requiring an extensive grid search to achieve high performance. By contrast, relatively simple diminishing stepsize rules are known to work well in subgradient-based approaches, both in theory and in practice. Thus, one may suspect that the ``double role'' of the $\eta$ parameter is what makes its online selection particularly challenging. Indeed, roughly speaking, the smaller $\eta$ the more ``local'' the direction is, taking into account only information computed close to the current iterate. However, towards termination, the algorithm should have accrued enough information that a ``global'' direction pointing towards the minimum of the model is likely preferable, which makes standard diminishing rules impractical.

\smallskip
\noindent
In order to improve both the performance and the robustness of existing approaches, we propose to delegate the selection of the convex multipliers and the choice of the stepsize to a novel \textit{objective-based} learning approach. Since our approach takes decisions dynamically during algorithm execution, Reinforcement Learning~\citep{sutton2018reinforcement} techniques could be, and in fact have been, successfully used~\citep{MAZYAVKINA2021,pmlr-v119-tang20a}. However, Reinforcement Learning is mainly useful when differentiating the execution is not possible, and it has many drawbacks such as being often sample-inefficient and suffering high-variance \citep{hassani2024sampleefficiencygeneralizationtransferinverse}. We rather propose to differentiate our algorithm via an unrolling approach~\citep{monga2020algorithmunrollinginterpretableefficient}, yielding a more direct optimization process with lower-variance gradient estimators by exploiting the internal solver dynamics \citep{scieur2023curseunrollingratedifferentiating,Kotary_2023}. A key challenge in this approach is to keep the structure of the \BM{} while differentiating through its execution. To address this problem, we use a heuristic approach whereby we replace these non-smooth components with smoother approximations based on a \NN{} model. Specifically, we propose a \NN{} that approximates the search direction using information from the current bundle and step size, replacing the $\eta$-strategy heuristic and the (dual) \MP{}~\eqref{DMP_nn}. This \NN{} integrates an approach inspired by the \BM{}, as outlined in Algorithm~\ref{pseudo:bm}, incorporating \textit{soft-updates} for the stabilization point, thereby enabling end-to-end differentiation using Automatic Differentiation~\citep{ADML} techniques.

\subsection{Literature Review}

Many recent studies have leveraged \ML{} to tackle Combinatorial Optimization problems~\citep{lodiTO2018}.
In this work, we focus on amortization techniques~\citep{semiAmort}, a research domain sometimes referred to as \textit{Learning-to-Optimize}~\citep{tianlong2022l2o, Li2016LearningTO,chen2022a}, and the associated methods are called \textit{semi-amortized optimization} methods.
In amortization techniques, the goal is to reduce the computational cost of \textit{iterative} optimization methods by learning a direct mapping from inputs, defining problem instances, to approximate solutions.
This often involves training an inference network to quickly produce a high-quality approximate solution that would otherwise require costly iterative optimization.
Many amortized frameworks using an iterative structure for predictions often require some \textit{unrolling technique}~\citep{monga2020algorithmunrollinginterpretableefficient} for back-propagation, involving reformulating the algorithm execution to enable the gradient computation in the backward pass.

Previous works designed to find the best parameters for an optimization algorithm have focused on finding the best step size for stochastic gradient descent~\citep{SinhaG71} in the context of Machine Learning, and are thus seen as instances of meta-learning
While early works considered simple evolutionary rules~\citep{bengio1990,bengio1995}, more recent ones are based on~\NNs{}~\citep{runarsson2000}.
Some works~\citep{Younger2001MetalearningWB,Hochreiter2001} jointly train with gradient descent both the networks for the learning problem and the one for the meta-learning task. 
However, these approaches still struggle to scale to modern architectures with tens of thousands of parameters.
To address this limitation, Andrychowicz et al.~\citep{NIPS2016_fb875828} introduce an LSTM-based optimizer that operates coordinate-wise on the network parameters.
This design enables different behavior on each coordinate by using separate activation functions for each parameter, while maintaining the optimizer's invariance to parameter order, as the same update rule is independently applied to each coordinate.

All of the above approaches fundamentally adopt the structure of gradient descent as the base for their learning framework. 
There is a growing interest in applying the \BM{} for optimization problems arising from the \ML{} community~\citep{le2007bundle}.
Variants of \PBMs~\citep{Krzysztof2006pbmal} are used to optimize complex, non-smooth loss functions often encountered in \ML{}~\citep{Alasdair2022sbmIN}. A specialized asynchronous bundle is presented for solving distributed learning problems~\citep{cederberg2025an}.
A further connection with the \BM{} can be found in \citep{ji2019learning}, where a meta-learning approach is proposed to coordinate learning in master-slave distributed systems.
This method aggregates the gradients for gradient descent, improving the scalability of distributed learning.

Recent research on \BMs{} for non-convex problems \citep{Oliveira2014LevelBM,de2016doubly}, opens the possibility to use the Bundle variants for \ML{} problems. \BMs{} have been applied to improve the training of SVMs~\citep{cortes1995support}, especially for large-scale problems~\citep{thorsten2009bundleSVM,teo2010bundleRM}.
Extensions of \BMs{} have also been proposed for multiclass SVMs, where \BMs{} can handle a larger number of constraints efficiently~\citep{thomas2008svm}.

The benefits of using a proximal operator when optimizing a convex non-smooth function have already been demonstrated~\citep{liu23e}. 
Also \citep{zhang2024metalearning} shows that advantages can be found in using proximality operators for meta-learning tasks, even if they did not consider the bundle method, but they use unrolling of the proximal gradient descent.

In \citep{li2024pdhgunrolledlearningtooptimizemethodlargescale}, a method is presented that approximates the resolution of a large-scale linear programming problem based on the associated Lagrangian function, providing iterative primal and dual updates using a neural network trained with unrolling techniques.

\subsection{Bundle Method}

The \BM{} (BM)~\citep{lemarechal1981bundle,lemarechal1995new} is an iterative method that solves problems of type
\begin{equation}
 \min \{ \, \phi(\vpi) \,:\, \vpi \in \Pi \, \}
 \label{generalconvexProblem}
\end{equation}
where $\phi : \mathbb{R}^m \rightarrow \mathbb{R}$ is a finite-valued convex non-smooth function and $\Pi$ is a ``simple'' convex set. The \LD{} is a relevant application of \eqref{generalconvexProblem}, discussed with more details in Section \ref{sec:lr}. In the rest of the paper we restrict ourselves to the case in which $\Pi = \sR^m$ or $\Pi = \sR^m_+ = \{ \vpi \in \sR^m \,:\, \vpi\geq 0 \}$, which correspond to the \LD{}~\citep{lrformilp}---that we use in our computational testing---when one is relaxing, respectively, equality or inequality constraints.
For the sake of presentation we only consider the case $\Pi = \sR^m$ here, but the approach can be simply extended, as discussed at the end of Section \ref{sec:lr}, to the case $\Pi = \sR^m_+$.

\smallskip
\noindent
Bundle methods iteratively optimize an approximation $\phi_{\beta}$ of the function $\phi$, constructed using the \emph{bundle} $\beta$ of first-order information collected so far, plus some regularization.
A common choice for $\phi_{\beta}$ is the \CP{} Model, i.e., the natural polyhedral upper approximation built using the very definition of subgradients.
For reasons to become apparent shortly, it is expedient to fix some reference point $\bar{\vpi}$; then, at a given iteration $t$ one has the set of previous iterates $\{ \vpi_i \}_{i=1}^t$ and can define the bundle $\beta_t$ as the set of pairs $\{ ( g_i , \alpha_i ) \}_{i=1}^t$, where $g_i \in \partial \phi(\vpi_i)$ is a subgradient obtained in the given iterate and $\alpha_i$ is the corresponding \emph{linearization error} $\alpha_i =  \phi(\bar{\vpi}) - \phi(\vpi_i) -\vg_i^{\top}( \bar{\vpi} - \vpi_i ) \geq 0$ (non-negativity being implicit into convexity of $\phi$).
Remarkably, the $\alpha_i$ depend on $\bar{\vpi}$, but we don't stress this notationally since $\bar{\vpi}$ will always be clear from the context.
In fact, $\beta_t$ may contain a strict subset of the thusly computed pairs (in case of bundle selection) and even pairs obtained differently (in case of bundle aggregation), but this is immaterial to our discussion, and we leave these details out for simplicity of presentation.
With this information, the \CP{} model~\citep{schrijver1980cutting} is
\begin{equation}\label{eq:bm}
    \phi_{\bar{\vpi},\beta_t}(\vd) = 
    \max\{ \, \vg_i^{\top} \vd + \alpha_i \,:
           \, ( g_i , \alpha_i ) \in \beta_t \, \} \;.
\end{equation}
It is easy to see that $\phi(\bar{\vpi} + \vd) - \phi(\bar{\vpi}) \geq \phi_{\bar{\vpi},\beta_t}(\vd)$, i.e., the \CPM{} is a polyhedral lower approximation of $\phi$ when translated w.r.t.~$(\bar{\vpi},\phi(\bar{\vpi}))$. Thus, however chosen $\bar{\vpi}$, the ``unstabilized'' \MP{}
\[
 \min \{ \, \phi_{\bar{\vpi},\beta_t}(\vd) \,:\,
            \bar{\vpi} + \vd \in \Pi \, \}
\]
would provide a lower bound on the optimal value of Problem~\eqref{generalconvexProblem}, but it is easy to see that such a problem can easily be unbounded below. Even if it is not, its optimal solution may easily be a very poor approximation to the true one, since $\phi_{\bar{\vpi},\beta_t}$ may be a poor model of $\phi$ ``far'' from the points in which the function has been computed. Because of this, one has to introduce some stabilization term that avoids wandering away from the region of space where $\phi$ is an appropriate model. The most common choice, corresponding to the \emph{proximal} \BM{}, is the Euclidean norm of the distance from the reference point $\bar{\vpi}$, which is therefore called the (current) \emph{stabilization centre}; this leads to the stabilized (primal) \MP{}
\begin{equation}
(MP_t) ~~
\min\{ \, \phi_{\bar{\vpi},\beta}(\vd) + \|\vd\|_2^2 \,/\, (2\eta) 
          \,:\, \bar{\vpi}+\vd \in \Pi \, \}
 \label{eq:PMP}
\end{equation}
whose optimal solution $\vd^*$ is used to choose the next iterate as $\vpi = \bar{\vpi} + \vd^*$. Solving~\eqref{eq:PMP} is equivalent to solving its dual, which for the case $\Pi = \sR^m$ reads
\begin{subequations} \label{DMP_nn}
\begin{align}
    \max & \textstyle
    \frac{\eta_{t}}{2} \big\| \sum_{(\vg_i,\alpha_i) \in \beta_t} \vg_i\theta_i \big\|_2^2+\sum_{(\vg_i,\alpha_i) \in \beta_t} \alpha_i\theta_i
    \\
    \text{s.t.} &
    \textstyle
    \sum_{(\vg_i,\alpha_i) \in \beta_t} \theta_i = 1 \;\;,\;\;
    \vtheta \geq 0 
\end{align} 
\end{subequations}
while the case $\Pi = \sR^m_+$ adds a new set of variables and simple constraints. The optimal solution $\vd^*_t$ of \eqref{eq:PMP} and the optimal solution $\vtheta^{(t)}$ of \eqref{DMP_nn} are linked by the relationship $\vd^{(t)} = \eta_{t} \vw_t$, where $\vw_t = \sum_{(\vg_i,\alpha_i) \in \beta_t} \vg_i \theta_i^{(t)}$ is the \emph{aggregated subgradient} computed out of the current bundle. This is relevant in that $\vw_t$ is a $\sigma_t$-subgradient of $\phi$ in $\bar{\vpi}$, where $\sigma_t = \sum_{(\vg_i,\alpha_i) \in \beta_t} \alpha_i \theta_i^{(t)}$ is the aggregated linearization error. For the case $\Pi = \sR^m_+$, $\vd^{(t)}$ is basically $\vw_t$ projected over the set of the active non-negativity constraints $\vpi \geq 0$; we refer to \citep[Chapter 4]{Lemarechal2001} for further details. Often, solving the \DMP{}~\eqref{DMP_nn} is more convenient than solving the primal due to the existence of ad-hoc approaches~\citep{frangioni1996solving}. The complete algorithm is described in Algorithm \ref{pseudo:bm}.

\begin{algorithm}
\caption{\BM{} pseudo-code}
\begin{algorithmic}[1]
\State \textbf{Choose} $\vpi_0,\eta_0,m$ 
\State $\bar{\vpi}_0 \gets \vpi_0$
\Comment The initial stabilization point coincides with the initial point
\State $(\vg_0,\valpha_0) \gets (\partial \phi(\vpi_0),0)$;
       $\beta_0 \gets \{(\vg_0,\valpha_0)\}$;
       $t \gets 0$
\While{\textbf{! stopping criteria}} \label{alg:sc}
    \State $(\vd^{(t)},\vtheta^{(t)}) \gets $ Solve $MP(\eta_{t},\beta_{t})$ \label{alg:DQP}
    \Comment Solve the \MP{}~\eqref{eq:PMP}--\eqref{DMP_nn}
    \State $\vpi_{t+1} \gets \bar{\vpi}_{t} + \eta_{t} \vd^{(t)}$
    \label{alg:tp}
    \Comment Compute the new trial point
    \State $(\vg_{t+1},\alpha_{t+1}) \gets (\partial \phi(\vpi_{t+1}),\phi(\bar{\vpi}_{t})-\phi(\vpi_{t+1})-\vg_{t+1}^{\top}(\bar{\vpi}_{t}-\vpi_{t+1}) )$ \label{alg:gle}
    
    \Comment Compute the new sub-gradient and linearization error
    \State $\beta_{t+1} \gets \beta_{t}\cup (\vg_{t+1},\valpha_{t+1})$
    \Comment Update the Bundle
  \label{alg:upd_B}
    \If{ $\phi(\vpi_{t+1})-\phi(\bar{\vpi}_{t}) > 
          m( \, \phi_{\beta_t}(\vpi_{t+1}) - \phi(\bar{\vpi}_{t}) \, )$ }
        \label{alg:cond_sp}
        \State $\bar{\vpi}_{t+1} \gets \vpi_{t+1}$;
               Update $\alpha_i\; \forall (\vg_i,\alpha_i) \in \beta_t$;
               Update $\eta_{t}$
               \label{alg:upd_sp} 
        \Comment Serious Step
     \Else
        \State $\bar{\vpi}_{t+1} \gets \bar{\vpi}_{t}$;
               Update $\eta_{t}$ \label{alg:nullStep}
    \Comment Null Step
    \EndIf
    \State $\beta_{t+1} \gets$ Remove outdated components $\beta_{t+1}$;
           $t \gets t+1$ \label{alg:roc}
\EndWhile
\end{algorithmic}
    \label{pseudo:bm}
\end{algorithm}

\smallskip
\noindent
The next trial point $\vpi_{t+1}$ is computed by making a step of $\eta_{t}$ along the direction $\vd_t$, which is either the (opposite of) the aggregated sugradient $\vw_t$ or easily obtained out of it. Computing $\phi(\vpi_{t+1})$ and its subgradient $\vg_{t+1}\in \partial\phi(\vpi_{t+1})$ enriches the bundle. The stabilization center is changed to $\vpi_{t+1}$, called a \textit{Serious Step}, if a ``substantial increase'' is obtained, measured with respect to the ascent $\phi_{\beta_t}(\vpi_{t+1})$ predicted by the model, with $m < 1$ a fixed parameter (usually ``small'', say $0.001$). Otherwise, i.e., if no (significant) progress is made, the stabilization center is not changed, which is referred to as a \textit{Null Step} (NS); yet, new information has been added to the bundle which will lead to a different direction $\vd_{t+1}$ in the next iteration. $\eta_{t}$ can (although it does not necessarily need to) also be adjusted, typically in different ways depending on the type of step taken. The stopping criterion in line~\ref{alg:sc} can be $\eta^* \| \vw_t \|_2^2 + \sigma_t \leq \epsilon \max\{0, \phi(\bar{\vpi}_t)\}$, where $\eta^*$ is ``large w.r.t.~$\eta_t$''; this corresponds to both $\vw_t \approx 0$ and $\sigma_t \approx 0$, proving that $\bar{\vpi}_t$ is an approximate minimum. At the end of each iteration, efficient implementations of the \BM{} may delete some old pairs from $\beta$ to lessen the \MP{} cost. A common strategy is removing components that are not used in the (dual) optimal solution over the last $K$ iterations (a \textit{quite} conservative strategy takes $K \sim 20$~\citep{frangioni2002generalized}). 

\section{Machine-Learning Based Bundle Method}\label{sec:objective_bn}

The approach presented in this work is heavily inspired by Algorithm \ref{pseudo:bm}, but with two major differences:
\begin{enumerate}
 \item we replace the solution of the \MP{}, which is a convex quadratic program that may become rather costly to solve as the size of the bundle grows, by the prediction of neural network and 
 \item we decouple the role of $\eta$ in the formation of $\vd$, via the weights $\vtheta$, from its role as a stepsize along it, by gain predicting its value with a neural network.
\end{enumerate}
To this effect, we need parametrize neural networks so that our two predictors from input instances and partial bundles to respectively convex multipliers $\vtheta$ and stepsize $\eta$ are accurate.
This parametrization is data-driven: we find the parameters that minimize a loss function, representing the value computed by the bundle network over a set of examples. 

\subsection{Learning Framework}

Let $\mathcal{D}$ be a dataset of instances $\iota$ and we denote by $\vpi_0\equiv\vpi_{0,\iota}$ the given initial point for the associated instance, possibly equal to zero. To mathematically define the learning problem, the number of iterations $T$ is fixed during the training process. Once the model is trained, the number of iterations is no longer a consideration at inference time~\citep{Li2016LearningTO}. In Section \ref{sec:lbn_bm} we describe in details the model $\vpi_T(\mW)\equiv \vpi_T(\mW;\vpi_0, \iota)$, where $\mW$ are the parameters. We may consider finding the parameters $\mW$ that, on average over instances drawn from $\mathcal{D}$, minimize the loss as defined by just the objective value $\mathcal{L}(\mW) = \phi(\vpi_T(\mW))$ after $T$ iterations:
\[
 \textstyle
 \min_{\mW} \mathbb{E}_{\iota\sim\mathcal{D}}[\mathcal{L}(\mW)] \;.
\]
However, Learning-to-Optimize approaches generally consider the contribution of each step of the whole trajectory, since it encourages optimizing the contributions of each step. Following this approach, we rather aim at maximizing the function:
\begin{equation}
 \textstyle
 \mathcal{L}_{\bm{\gamma}}(\mW) = \sum_{t=1}^T \gamma_t\phi(\vpi_t(\mW))
\end{equation}
where $\vpi_t(\mW)$ denotes the trial point at iteration $t$, which is computed by a \NN{}. We use $\gamma_t=\gamma^{T-t}$, where $0< \gamma\leq 1$ is a hyperparameter. Hence, the learning problem can be formulated as
\begin{equation}
 \textstyle
 \min_{\mW} \mathbb{E}_{\iota\sim\mathcal{D}}[\mathcal{L}_{\bm{\gamma}}(\mW)].    
\end{equation}
In the next sections, we will specify how to differentiate the loss function using the chain rule and the structure of the model $\vpi_{t}(\mW)$ used as \ML{} surrogate of the \MP{} in the \BM{}.

\subsection{Loss Function}

We want to compute a subgradient $\partial_{\mW}\mathcal{L}_{\vgamma}$ of the loss with respect to the network parameters. 
$\mathcal{L}_{\bm{\gamma}}(\mW)$ is a linear combination of the function values $\phi$ evaluated at different points along the trajectory: consequently, its subgradient with respect to $\mW$ can be written as the weighted sum of subgradients of the function $\phi$ computed in different points. Thus we only need to compute subgradients $\partial_{\mW}\phi(\vpi_t(\mW))$ for each $t \leq T$.

\noindent
We apply the chain rule to compute $\partial_{\mW}\phi(\vpi_t(\mW)) $ as $\partial_{\vpi_t}\phi(\vpi_t(\mW))\cdot\partial_{\mW}\vpi_t(\mW)$. Note that $\partial_{\vpi}\phi(\vpi_t(\mW))$ is the subgradient of $\phi$ at $\vpi_t(\mW)$, which exists since $\phi$ is convex and in fact corresponds to the $\vg_t$ used in the standard version of the \BM{}. As we construct $\vpi_t(\mW)$ to be differentiable, $\partial_{\mW}\vpi_t(\mW)$ corresponds to the gradient $\nabla_{\mW}\vpi_t(\mW)$ that can be obtained using automatic differentiation techniques. The algorithm structure for this computation is discussed in detail in the next section. 

Putting all together, we obtain the following gradient approximation
\[
\textstyle
\partial_{\mW}\mathcal{L}(\mW)=\sum^T_{t=1}\gamma^{T-t}\vg_t\nabla_{\mW}\vpi_t(\mW)
\]

\subsection{Unrolled Model}\label{sec:lbn_bm}

In this section, we provide further details on the structure of \ML{} model used to compute the trial point $\vpi_t(\mW)$ at iteration $t$. We keep the structure of the \BM{} mostly unchanged, modifying only some operations in order to make them differentiable: the choice of the step size, the solution of Problem \ref{DMP_nn}, and the update of the stabilization center. The result is presented in Algorithm \ref{alg:unrolled_bm}: since it only employs differentiable operations, we can then use standard automatic differentiation tools.

\noindent
At each iteration, we need to provide a feature vector to the neural network. We do not back-propagate through the feature extraction, which is made using hand-crafted features described in detail in Appendix~\ref{sec:hyper_parameters_and_implementation}. These features are fed as input to a \NN{}, which outputs a vector $\vtheta^{(t)}\in \Delta^{t}=\big\{ \vtheta \in [0,1]^{t+1} \; | \; \sum_{i=1}^{t+1} \theta_i = 1\big\}$ in the $t-$dimensional simplex and a step size $\eta_t$. More details on the network architecture can be found in Section~\ref{sec:bn_numRes}.

In the standard \BM{}, the stabilization center is updated if the condition $\phi(\vpi_{t+1}) > \phi(\bar{\vpi}_{t})$ is satisfied. However, this is a piecewise-constant operation, and therefore with zero derivative almost everywhere. Denoting by $\vr=\text{arg}\min((\phi(\vpi_{t+1}),\phi(\bar{\vpi}_{t})))$ the one-hot vector encoding equal to $(1,0)$ if $\phi(\vpi_{t+1})\geq\phi(\bar{\vpi}_{t})$ and $(0,1)$ otherwise, the update of the stabilization point can be rewritten as
\[
\vr^{\top}(\vpi_{t+1}, \bar{\vpi}_{t})=r_1 \vpi_{t+1}+r_2 \bar{\vpi}_{t}.
\]
This corresponds to choosing a vertex of the two-dimensional simplex. A common strategy in \ML{} to obtain smoother approximations of the $\text{arg}\max$ operator is to replace it with a \textit{soft} version, obtained with the $\text{softmin}$~\cite[pp. 180-184]{goodfellow2016deep}. This leads to considering
\[
\vr=\mathbf{softmin}(\phi(\vpi_{t+1}),\phi(\bar{\vpi}_{t}))
\]
and the updates can be performed as previously, i.e., choosing a possibly different convex combination of the two vectors $\vpi_{t+1}$ and $\vpi_{t}$, which will not necessarily be a vertex of the two-dimensional simplex, but a point inside it.

It is important to notice that we do not differentiate through the operation computing the gradient in line \ref{line:grad_comp}, as obtaining exact derivatives requires the function $\phi$ to be twice differentiable. It is possible to approximate this contribution, but this is not explicitly considered in this work; this strategy has been commonly pursued in other works, for example \citep{NIPS2016_fb875828}. In other words, the bundle information fed to the neural network is considered as a constant, even though it is produced by a gradient computation. This is similar to the way recurrent language models do not backpropagate through the previously generated words, but rather through the previous latent state~\citep{bahdanau-2015-neural-machin}.

\begin{algorithm}[h]
\caption{Bundle Network pseudo-code}
\begin{algorithmic}[1]
\State \textbf{Choose} $\vpi_0,T$
\State $\bar{\vpi}_0 \gets \vpi_0$ \Comment{Initialize stabilization point}
\State $(\vg_0, \valpha_0,v_0) \gets (\partial \phi(\vpi_0), 0,\phi(\vpi_0))$
\State $\beta_0 \gets \{(\vg_0, \valpha_0)\}$
\For{$t=1,\cdots,T$}
          \State $ \varphi_{t} = $ features\_extraction $(\beta_{t})$
          \Comment{We do not back-propagate through this operation}
        \State $\eta_t, \bm{\delta}^{(t)} \gets nn(\varphi_t)$  \Comment{Output of a \NN{}}
         \State $\vtheta^{(t)} \gets \psi(\bm{\delta}^{(t)})$  \Comment{Approximate Solution of the DMP~\ref{DMP_nn}}
        \State $\vw^{(t)} \gets \sum_{i=0}^{|\beta_t|}\vg_i\theta_i^{(t)}$ \Comment{New trial direction}
        \State $\vpi_{t+1} \gets \sigma(\bar{\vpi}_{t} + \eta_{t} \vw^{(t)})$ \Comment{Compute new trial point}\label{line:generalization}
        \State $(\vg_{t+1}, \valpha_{t+1},v_{t+1}) \gets (\partial \phi(\vpi_{t+1}), \phi(\bar{\vpi}_{t}) - \phi(\vpi_{t+1}) - \vg_{t+1}^{\top}(\bar{\vpi}_{t} - \vpi_{t+1} ),\phi(\vpi_{t+1}))$ \label{line:grad_comp}\Comment{Gradient and linearization error}
        \State $\beta_{t+1} \gets \beta_{t} \cup (\vg_{t+1}, \valpha_{t+1},v_{t+1})$ \Comment{Update Bundle}
         \State $\bar{\vpi}_{t+1} \gets \mathbf{softmin}(\phi(\vpi_{t+1}),\phi(\bar{\vpi}_{t}))\odot (\vpi_{t+1},\bar{\vpi}_{t})$ 
         \Comment{Smooth updates of the stabilization point}
                 \State Update $\alpha_i\; \forall i =0,\cdots,|\beta_t|$
\Comment{Update the linearization errors}
    \EndFor
\end{algorithmic}
\caption{Machine learning-based version of the \BM{}~\ref{pseudo:bm} that substitutes the resolution of the \MP{} and the non-continuous operations with a smoother version based on \NNs{}.}
\label{alg:unrolled_bm}
\end{algorithm}

\subsubsection{Step and Direction Architecture}\label{sec:overall_architecture_bn}

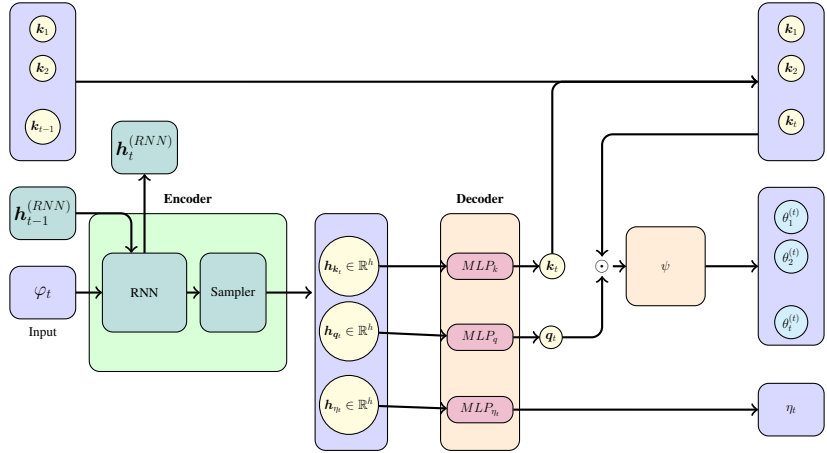
\begin{figure}[h]
    \centering
    \begin{tikzpicture}[node distance=3cm, scale=0.35, transform shape]

\tikzset{
    c_data/.style={rectangle, draw=black, fill=blue!15, minimum width=2.5cm, minimum height=2cm, align=center},
    r_flow/.style={rectangle, draw=black, fill=teal!25, minimum height=3cm, minimum width=2.5cm, align=center},
    rr_flow/.style={rectangle, draw=black, fill=orange!15, minimum height=4cm, minimum width=3cm, align=center},
    r_mlp/.style={rectangle, draw=black, fill=purple!25, minimum width=2.5cm, minimum height=1cm, align=center},
    z_tiny/.style={circle, draw=black, fill=yellow!15, minimum size=0.5cm},
    pi_tiny/.style={circle, draw=black, fill=cyan!15, minimum size=0.75cm},
    arrow/.style={->, thick},
}

\node (data) [c_data,rounded corners] {\huge  $\varphi_{t}$};
\node (data_labels) [below=0.2cm of data.south,rounded corners] {\Large Input};
   \node (Ks) [c_data, above=4cm of data,minimum height=6cm,rounded corners] {};
    
\foreach \i/\y in {1/0.5, 2/2} {
    \node (k_\i) [z_tiny, below=\y of Ks.north,rounded corners] {\Large$\vk_{\i}$};
}
    \node (k_tm1) [z_tiny, below=4.05 of Ks.north,rounded corners] {\Large$\vk_{t-1}$};

\node (encoder) [rectangle, minimum height=6cm, minimum width=7.5cm, right=0.5cm of data, draw=black, fill=green!15,rounded corners] {};
\node (encoder_labels) [above=0.25cm of encoder.north,rounded corners] {\Large \textbf{Encoder}};

\node (NN1) [r_flow, right=1cm of data, minimum width=3.2cm,rounded corners] {\Large RNN};

\node (Sampler) [r_flow, right=0.5cm of NN1,rounded corners] {\Large Sampler};

\node (hrnntm1) [r_flow, above = 1cm of data,rounded corners,minimum height=2cm] {\huge  $\vh_{t-1}^{(RNN)}$};
\node (hrnnt) [r_flow, above = 3cm of NN1,rounded corners,minimum height=2cm] {\huge  $\vh_{t}^{(RNN)}$};


\node (kq_box) [c_data,right=1.85cm of Sampler,minimum height=9cm,minimum width=2.75cm,rounded corners,yshift=-1.5cm] { };
    \node (kt) [z_tiny,right=2cm of Sampler,yshift=1cm] {\Large $\vh_{\vk_t}\in \sR^h$};
    \node (qt) [z_tiny, below=0.2cm of kt] {\Large $\vh_{\vq_t}\in \sR^h$};
    \node (etat) [z_tiny, below=3cm of kt] {\Large $\vh_{\eta_{t}}\in \sR^h$};

\node (mlps) [rr_flow, right=2cm of kq_box,minimum height=9cm,rounded corners] {};

\foreach \i/\y in {k/1.5, q/4.2} {
    \node (mlp_\i) [r_mlp, below=\y of mlps.north,rounded corners] {\Large$MLP_{\i}$};
}
    \node (mlp_eta) [r_mlp, below=6.95 of mlps.north,rounded corners] {\Large$MLP_{\eta_t}$};

\node (encoder_labels) [above=0.25cm of mlps.north,rounded corners] {\Large \textbf{Decoder}};

    \node (k_t) [z_tiny, right=1cm of mlp_k,rounded corners] {\Large$\vk_{t}$};

 \node (q_t) [z_tiny, right=1cm of mlp_q,rounded corners] {\Large$\vq_{t}$};

\node (dp) [right=1cm of k_t,rounded corners] {\huge$\odot$};

\node (sm) [rr_flow,right=0.5cm of dp,minimum height=3cm,rounded corners] {\Large $\psi$};

\node (post_sm) [c_data, right=2cm of sm,minimum height=6cm,rounded corners] {};

\foreach \i/\y in {1/0.5, 2/2, t/4.5} {
    \node (theta_\i) [pi_tiny, below=\y of post_sm.north] {\Large$\theta_\i^{(t)}$};
}

\node (t) [c_data, rounded corners] at (mlp_eta.east -| post_sm.south) {\Large $\eta_{t}$};

 \node (Ks_new) [c_data,minimum height=6cm,rounded corners]  at (Ks.east -| post_sm.north) {};
   
\foreach \i/\y in {1/0.5, 2/2,t/4.05} {
    \node (kn_\i) [z_tiny, below=\y of Ks_new.north,rounded corners] {\Large$\vk_{\i}$};
}

    \draw [arrow,rounded corners] (mlp_q) --  (q_t) ;
    \draw [arrow,,rounded corners] (mlp_k) --  (k_t)  ;

\draw [arrow] (Sampler) -- ++(30mm,0mm);
\draw [arrow] (sm) -- (post_sm);
\draw [arrow]  (qt) --  (mlp_q);
\draw [arrow] (kt.east) |- (mlp_k);
\draw [arrow,rounded corners]  (q_t) -| (dp);

\draw [arrow] (dp) -- (sm);

\draw [arrow] (data) -- (NN1);
\draw [arrow] (NN1) -- (Sampler);

\draw [arrow,rounded corners] (Ks_new.west)+(0,-2cm) -|  (dp);

\draw [arrow,rounded corners] (Ks) -- (Ks_new);

\draw [arrow,rounded corners] (k_t.north) |- (Ks_new.west);


\draw[arrow,rounded corners] (hrnntm1) -| ($(NN1.north)+(-0.5cm,0)$);
\draw[arrow,rounded corners]  (NN1.north) -- (hrnnt) ;

\draw [arrow,rounded corners] (etat) -- (mlp_eta);

\draw [arrow,rounded corners] (mlp_eta) --  (t);


\end{tikzpicture}

    \caption{Schematic representation of the architecture used to predict the step-size $\eta_t$ and the weights $\vtheta^{(t)}$ for the convex combination of the gradients in the bundle at the current iteration.}
    \label{fig:BundleNetwork}
\end{figure}

This section provides further details on the architecture of the \ML{}-model, schematically presented in Figure~\ref{fig:BundleNetwork}, which is used to predict $\eta_t$ and $\vtheta^{(t)}$.

\noindent
Given the features $\varphi_{t}\in\sR^d$ at the current bundle $\beta_{t}$, a \RNN{}~\citep{chung2014empirical,cho2014learning,kyunghyun2014}, specifically an \LSTM{} (LSTM~\citep{lstmHoc}) layer, maps this vector to a hidden space by predicting six vectors of size $h$. We denote these vectors as $\bm{\mu}_{{q_t}},\bm{\sigma}_{{q_t}},\bm{\mu}_{{k_t}},\bm{\sigma}_{{k_t}},\bm{\mu}_{{\eta_t}},\bm{\sigma}_{{\eta_t}}\in  \sR^h$, and they represent the mean and the variance used to extract, by sampling through a Gaussian distribution, the hidden representation of the keys and the queries used to represent the current iterations as in attention mechanisms~\citep{vaswani2017attention}, as well as the hidden representation for $\eta_t$. More precisely, to obtain the hidden representations $\vh_{q_t},\vh_{k_t},\vh_{\eta_t}\in \sR^h$, we consider a Gaussian distribution, so we can approximate the expectation by sampling and use the reparameterization trick~\citep{Kingma2014, NIPS2015_de03beff} to perform standard backpropagation.
The size of the hidden representations and other hyperparameters of the network will be discussed in Appendix~\ref{sec:hyper_parameters_and_implementation}.

The sampled representations are then fed into three independent \MLPs{} that predict a scalar $\eta_t$ and two vectors $\vq_t$ and $\vk_t$.
The vectors $\vq_t$ and $\vk_t$ are, respectively, the query and the key of the current iteration.
The vector $\vk_t$ is stored in memory for use in subsequent iterations.
Using $\vq_t$ and all the vectors $\{\vk_i\}_{i=0}^{|\beta_t|}$, we can compute a scalar for each component in the bundle using an \textit{Attention Layer} defined as follows:
$$
\bm{\delta}^{(t)}=(\vk_i^{\top}\vq_t)_{i=0}^{|\beta_t|} \in \sR^t.
$$
To obtain the approximate solution of the \DMP{}~\ref{DMP_nn}, we pass this vector to a normalization function $\psi$ taking as input a vector and returning one of the same dimension, but with all components between zero and one, and summing to one
$$
\vtheta^{(t)}=\psi(\bm{\delta}^{(t)})\in \Delta^t.
$$
This yields a vector $\vtheta^{(t)}$ whose components lie in $[0,1]$ and sum to one, as the solution of the Problem \eqref{DMP_nn} considered in the \textit{classic} bundle. Some examples of $\psi$ are the softmax or sparsemax~\citep{MartinsA16} functions.

\section{Evaluation}\label{sec:bn_numRes}

In this section, we present the main numerical results of our work. We test it only on the example case of the \LR{}, of which we then provide further details.

\subsection{Lagrangian Relaxation}\label{sec:lr}

Let $\iota$ be a Mixed Integer Linear Program~\citep[Chap.~8]{confortiIntegerProgramming2014} of the form:
\begin{subequations}
\label{eqn:lr_formulation}
\begin{align}
	\label{eq:MILP} (P) \qquad \min & \; \vw^{\top}\vx \\
	\label{ctn:ax<=b}  & \mA\vx = \vb \\
	\label{ctn:cx<=d}  &  \mC\vx \ge \vd
    \;\;,\;\; \vx \in \mathbb{R}_+^{n_1}\times\mathbb{N}^{n_2} 
\end{align}
\end{subequations}
The \textit{\LP{}} is obtained by dualizing a family of constraints. Here we relax the constraints defined by the inequalities~\eqref{ctn:ax<=b}, and penalize their violation with \LMs{} (LMs) $\vpi \in \sR^{m}$:
\begin{align*}
	\big(LR (\vpi)\big) \qquad \min \; 
    & \vw^{\top}\vx+\bm{\pi}^{\top} (\vb-\mA\vx)  \\
	& \mC\vx \geq \vd{} \;\;,\;\;
	   \vx{}\in{} \mathbb{R}_+^{n_1}\times\mathbb{N}^{n_2}
\end{align*}
Weak Lagrangian duality ensures that $LR(\vpi)$ provides a lower bound for $P$. The goal of the \LDP{} is to determine the optimal bound.
\begin{equation}\label{eq:dual}
	(LD)  \qquad \max LR(\bm{\pi}).
\end{equation}
\LR{} is particularly valuable when removing certain constraints results in a \LP{} $LR(\vpi)$ that can be solved efficiently; it is therefore crucial to carefully select the set of constraints to dualize. Indeed, this choice involves a trade-off between the speed and the quality of the resulting bound, as will be discussed later, due to the \textit{No Free Lunch Theorem}~\citep{geoffrion_lagrangean_1974} that proves that, in order to achieve better bounds than the \CR{}, the \LP{} should be ``harder'' than a \Lp{}. Actually, the \LDP{} \ref{eq:dual} is a concave maximization problem, but it is trivial to reformulate it as a convex minimization problem~\citep{rockafellasCA}. Relaxing inequality constraints $\mA\vx \ge \vb$ instead of equality ones leads to nonnegative multipliers $\vpi\in\sR_+^m$. To ensure that we predict only nonnegative multipliers we substitute line \ref{line:generalization} of Algorithm \ref{alg:unrolled_bm} with
$$
\vpi_{t+1}=\sigma(\bar{\vpi}_t + \eta_t \vw^{(t)})
$$
where $\sigma$ is a component-wise non-negative activation function such as ReLU.

\subsection{Experiments}

In this section, we compare our approaches to classical iterative methods for solving the \LD{}, specifically those based on subgradient information. In particular, we evaluate the performance of the proposed approach against both simple subgradient-based methods and the classical \BM{}, using different heuristic $\eta$-strategies.

We consider ADAM~\citep{kingma2014adam} and Gradient Ascent using a simple strategy for regularizing the step size: if we pass more than two iterations without improving the objective value, we divide the step size by two. These two baselines are well-established in the context of convex minimization and serve as relevant points of comparison. They are of particular interest because they are typically faster than the \BM{}, as Descent relies only on the subgradient at the latter visited point, meanwhile, ADAM relies on the subgradient of the last inserted point and the previous search direction.

We further compare our approach with \textit{classic} variants of the \BM{} considering the four long-term $\eta$-strategies presented in Appendix~\ref {app:eta_strat}, inspired by the ones used in the state-of-the-art bundle-based \textit{BundleSolver} component of \textsc{SMS++}~\citep{smspp}. With the non-constant ones, we also consider the middle-term $\eta$-strategy and the short-term $\eta-$strategy proposing increasing as $\eta_{t+1}=1.1\cdot\eta_t$ and decreasing as $\eta_{t+1}=\eta_{t}\cdot 0.9$. For all these six \textit{classic approaches}, we perform a grid search to choose the best initial parameter $\eta_0$ considering different values. We evaluate the initialization of the $\eta_0$ parameter using the values $10^k$ for $k \in \{4, 3, 2, 1, 0, -1\}$. For each fixed maximum number of iterations ($10$, $25$, $50$, and $100$), we select the value of $\eta_0$ that achieves the lowest average GAP over the test set. The results of this grid search are summarized in Table~\ref{tab:grid_search}.

The objective values of $\phi$ are too large for a learning method, so we rescale the objective $\phi$ by dividing the function value using the norm $\sqrt{\vg^{\top}_0\vg_0}$ of the gradient $\vg_0$ in the starting point $\vpi_0$, always taken equal to the zero vector. Even if the objective value changes, the optimum of the original function $\phi$ and of its scaled version are the same.

We consider the \LD{} of the Multi-Commodity Network Design (MC) presented in Appendix \ref{app:MC} and the \LD{} of the Generalized Assignment (GA) presented in Appendix \ref{app:GA}. For MC, we consider a \LR{} in which we dualize equality constraints, while for GA we consider a \LR{} in which we dualize inequality constraints, applying the modifications previously discussed. For the MC and the GA, we use the datasets considered in ~\citep{demelas24}.

For each instance $\iota$ on which we evaluate our model, the \LR{} is solved using SMS++. The obtained \LMs{}, denoted by $\vpi^*_{\iota}$, are used for computing the GAP. We enhance the fact that we do not need these multipliers to perform the training, but we use them to provide a metric for our approach. We use the percentage gap as a metric to evaluate the quality of the bounds computed by the different systems, averaged over a dataset of instances \(\mathcal{I}\). We measure the quality of each approach by means of the percentage GAP
\[
 100 \times \frac{1}{| \mathcal{I} |}\sum_{\iota\in\mathcal{I}}\frac{LR(\bm{\pi}^*_{\iota})- B_{\iota} }{LR(\bm{\pi}^*_\iota)} \;,
\]
where $B_\iota$ is the returned upper bound for instance $\iota$.

\subsection{Bound Accuracy}\label{sec_main_table}
Table~\ref{tab:mainTable_bn} reports the performance of different systems on our datasets.
For the heuristic strategies that did not use \ML{}, we consider a grid search for the initialization of the step-size/regularization parameter.
We consider values in between $10^{4}$ and $10^{-1}$.
For each fixed number of iterations ($10,\; 25,\; 50$ and $100$), we choose the best $\eta_0$ concerning the values in the grid.
From Table~\ref{tab:mainTable_bn} we can see that Bundle Network outperforms other techniques in terms of GAP for a fixed number of iterations in different datasets.

\begin{table}[h!]
\centering
    \caption{New Comparison of our approach and different baselines. The baselines' starting regularization parameter/step size is chosen using a grid search. Meanwhile, our approach does not need hyperparameter tuning. Ours is trained only on $10$ iterations.}
    \begin{tabular}{@{}ll rr rr rr rr@{}}
    \toprule
\multirow{2}{*}{} & Methods 
& \multicolumn{2}{c}{10 iter.} 
& \multicolumn{2}{c}{25 iter.} 
& \multicolumn{2}{c}{50 iter.} 
& \multicolumn{2}{c}{100 iter.} \\
\midrule
& &  GAP  & time  
& GAP  & time   
& GAP  & time   
& GAP  & time  \\
\midrule
\multirow{8}{*}{\rotatebox[origin=c]{90}{\textsc{MC-SML-40}}} 
  & Bundle h. & 20.18 & 0.044 & 6.68 & 0.117 & 2.84 & 0.269 & 1.17 & 0.767\\
 & Bundle b. & 17.68 & 0.054 & 6.09 & 0.145 & 1.84 & 0.327 & 0.41 & 0.898\\
 & Bundle s. & 17.58 & 0.044 & 6.10 & 0.116 & 1.82 & 0.271 & 0.36 & 0.778\\
 & Bundle c. & 11.96 & 0.048 & 4.16 & 0.128 & 1.30 & 0.293 & 0.31 & 0.816\\
 & Descent & 44.66 & \textbf{0.008} & 14.11 & \textbf{0.021} & 4.70 & \textbf{0.044} & 2.24 & \textbf{0.093} \\
 & Adam & 48.87 & 0.011 & 12.84 & 0.027 & 3.04 & 0.055 & 1.89 & 0.109 \\
& Bundle n. & \textbf{9.20} & 0.105 & \textbf{1.75} & 0.185 & \textbf{0.53} & 0.321 & \textbf{0.17} & 0.592\\
\midrule
\multirow{8}{*}{\rotatebox[origin=c]{90}{\textsc{MC-SML-Var}}} 
 & Bundle h. & 18.27 & 0.094 & 9.55 & 0.254 & 4.08 & 0.564 & 1.96 & 1.433\\
 & Bundle b. & 15.83 & 0.124 & 7.39 & 0.334 & 2.88 & 0.736 & 0.74 & 1.749\\
 & Bundle s. & 15.82 & 0.094 & 6.92 & 0.256 & 2.58 & 0.576 & 0.67 & 1.421\\
 & Bundle c. & {12.79} & 0.108 & 4.33 & 0.29 & 1.66 & 0.635 & 0.58 & 1.503\\
 & Descent & 62.43 & \textbf{0.025} & 28.64 & \textbf{0.066} & 11.14 & \textbf{0.142} & 4.00 & \textbf{0.303} \\
 & Adam & 51.48 & 0.034 & 15.42 & 0.086 & 4.74 & 0.17 & 2.75 & 0.341 \\
& Bundle n. & \textbf{12.60} & 0.138 & \textbf{3.28} & 0.269 & \textbf{1.24} & 0.496 & \textbf{0.50} & 0.952\\
\midrule
\multirow{8}{*}{\rotatebox[origin=c]{90}{\textsc{MCND-Big-40}}} 
 & Bundle h. & 20.66 & 0.07 & 8.54 & 0.185 & 5.32 & 0.410 & 2.88 & 0.986\\
 & Bundle b. & 20.67 & 0.084 & 7.37 & 0.222 & 3.17 & 0.485 & 1.24 & 1.209\\
 & Bundle s. & 20.67 & 0.07 & 7.37 & 0.184 & 3.17 & 0.407 & 1.19 & 1.042\\
 & Bundle c. & 19.34 & 0.076 & 5.64 & 0.203 & 2.08 & 0.447 & 0.78 & 1.110\\
 & Descent & 28.53 & \textbf{0.023} & 14.12 & \textbf{0.058} & 9.91 & \textbf{0.117} & 8.21 & \textbf{0.241} \\
 & Adam & 24.89 & 0.025 & 9.38 & 0.064 & 7.30 & 0.129 & 6.90 & 0.256 \\
& Bundle n. & \textbf{12.70} & 0.130 & \textbf{3.41} & 0.243 & \textbf{1.26} & 0.436 & \textbf{0.47} & 0.832\\
  \midrule
\multirow{8}{*}{\rotatebox[origin=c]{90}{\textsc{MC-Big-Var}}} 
 & Bundle h. & 26.23 & 0.095 & 9.99 & 0.248 & 4.41 & 0.550 & 2.23 & 1.386\\
 & Bundle b. & 21.82 & 0.123 & 10.01 & 0.324 & 3.46 & 0.695 & 1.08 & 1.652\\
 & Bundle s. & 21.83 & 0.094 & 9.53 & 0.253 & 3.43 & 0.550 & 1.02 & 1.353\\
 & Bundle c. & \textbf{17.98} & 0.112 & 6.34 & 0.289 & 2.66 & 0.627 & 0.93 & 1.477\\
 & Descent & 54.11 & \textbf{0.029} & 24.69 & \textbf{0.074} & 10.58 & \textbf{0.154} & 4.85 & \textbf{0.323} \\
 & Adam & 44.21 & 0.034 & 13.75 & 0.086 & 5.42 & 0.170 & 3.84 & 0.345 \\ 
& Bundle n. & 20.10 & 0.139 & \textbf{5.11} & 0.271 & \textbf{2.01} & 0.495 & \textbf{0.70} & 0.953\\
\midrule
\multirow{8}{*}{\rotatebox[origin=c]{90}{\textsc{GA-10-100}} }
 & Bundle h. & 0.2333 & 0.070 & 0.028 & 0.241 & 0.0016 & 0.657 & \textbf{0.0006} & 1.596\\
 & Bundle b. & 0.2333 & 0.073 & 0.028 & 0.249 & 0.0016 & 0.684 & \textbf{0.0006} & 1.655\\
 & Bundle s. & 0.3055 & 0.065 & \textbf{0.0156} & 0.211 & 0.0026 & 0.553 & 0.0024 & 1.283\\
 & Bundle c. & {0.1893} & 0.071 & 0.0157 & 0.235 & \textbf{0.0014} & 0.601 & 0.0011 & 1.383\\
 & Descent & 0.8799 & 0.013 & 0.2091 & 0.032 & 0.0644 & 0.063 & 0.0364 & 0.127\\
 & Adam & 0.7234 & \textbf{0.012} & 0.1843 & \textbf{0.029} & 0.0183 & \textbf{0.059} & 0.0048 & \textbf{0.119}\\
& Bundle n. & \textbf{0.1484} & 0.104 & {0.0228} & 0.177 & {0.0047} & 0.304 & {0.0009} & 0.551\\
 \midrule
\multirow{8}{*}{\rotatebox[origin=c]{90}{\textsc{GA-20-400}}}
 & Bundle hard & 0.157 & 0.445 & 0.0343 & 1.548 & 0.0088 & 4.528 & 0.0021 & 15.099\\
 & Bundle balancing & 0.157 & 0.463 & 0.0343 & 1.568 & 0.0088 & 4.539 & 0.0021 & 15.158\\
 & Bundle soft & 0.3542 & 0.431 & 0.0284 & 1.417 & \textbf{0.0048} & 4.256 & 0.0028 & 14.145\\
 & Bundle constant & 0.1129 & 0.433 & 0.0193 & 1.459 & 0.0052 & 4.264 & 0.0020 & 14.459\\
 & Descent & 0.5815 & 0.205 & 0.2971 & 0.5 & 0.0753 & 0.98 & 0.0255 & 1.937\\
 & Adam & 0.6598 & \textbf{0.177} & 0.0831 & \textbf{0.444} & 0.016 & \textbf{0.887} & 0.0090 & \textbf{1.775}\\
& Bundle n. & \textbf{0.1018} & 0.228 & \textbf{0.0190} & 0.498 & {0.0052} & 0.943 & \textbf{0.0014} & 1.837\\
\bottomrule
\end{tabular}
    \label{tab:mainTable_bn}
\end{table}

\begin{figure}
    \centering
   \includegraphics[width=0.45\textwidth]{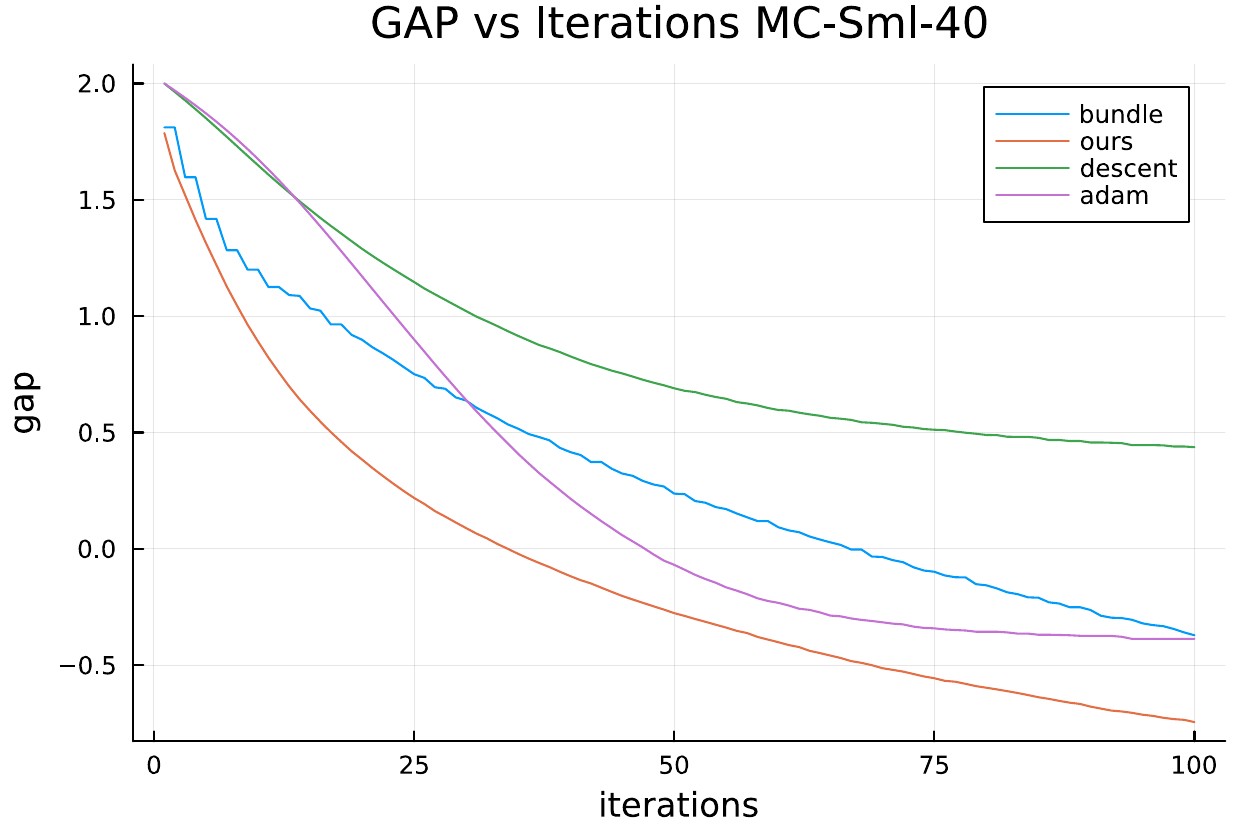}
   \includegraphics[width=0.45\textwidth]{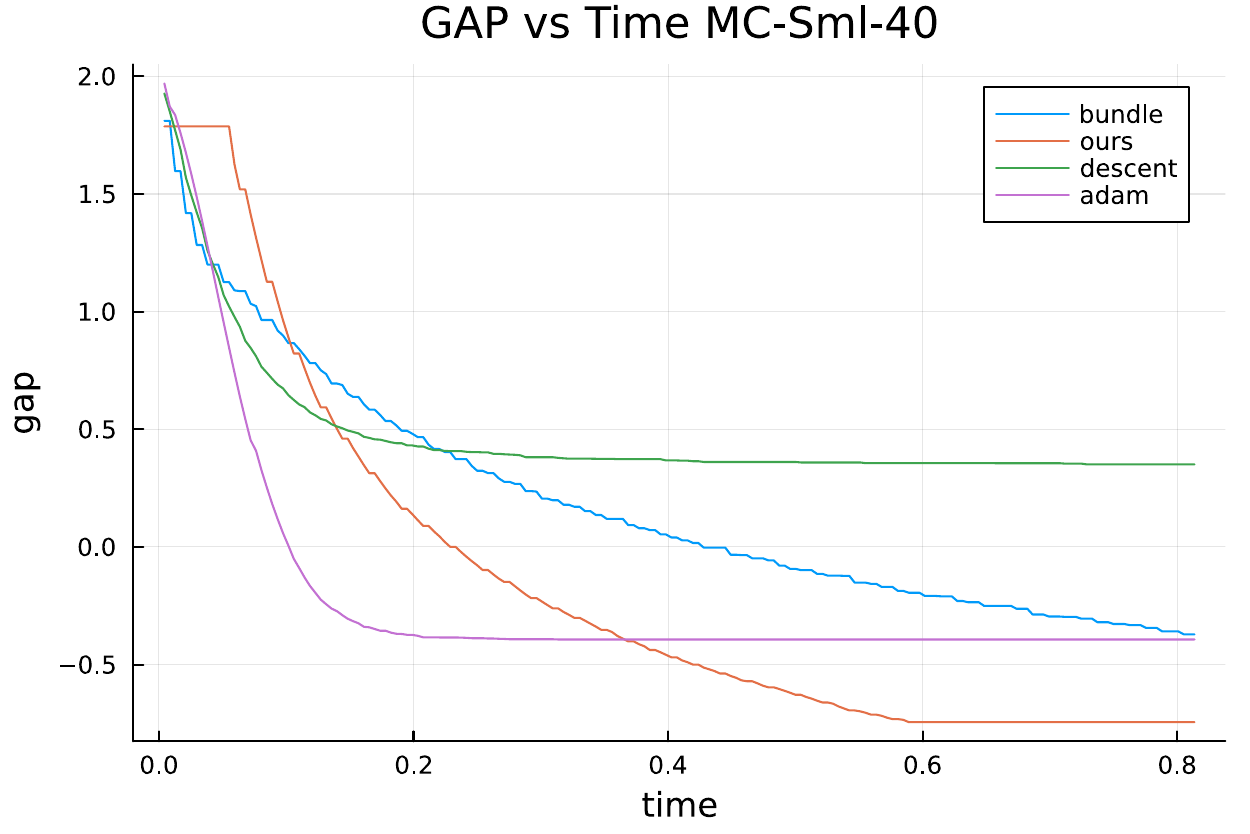}
    \caption{Comparison plots of the GAP (in logarithmic scale) during the execution of our trained \NN{} and different baselines in terms of iterations (in the left) and in terms of time (in the right) for \textsc{MC-Sml-40} dataset.}
    \label{fig:execution_MC_Sml_40}
\end{figure}

\begin{figure}
    \centering    
   \includegraphics[width=0.45\textwidth]{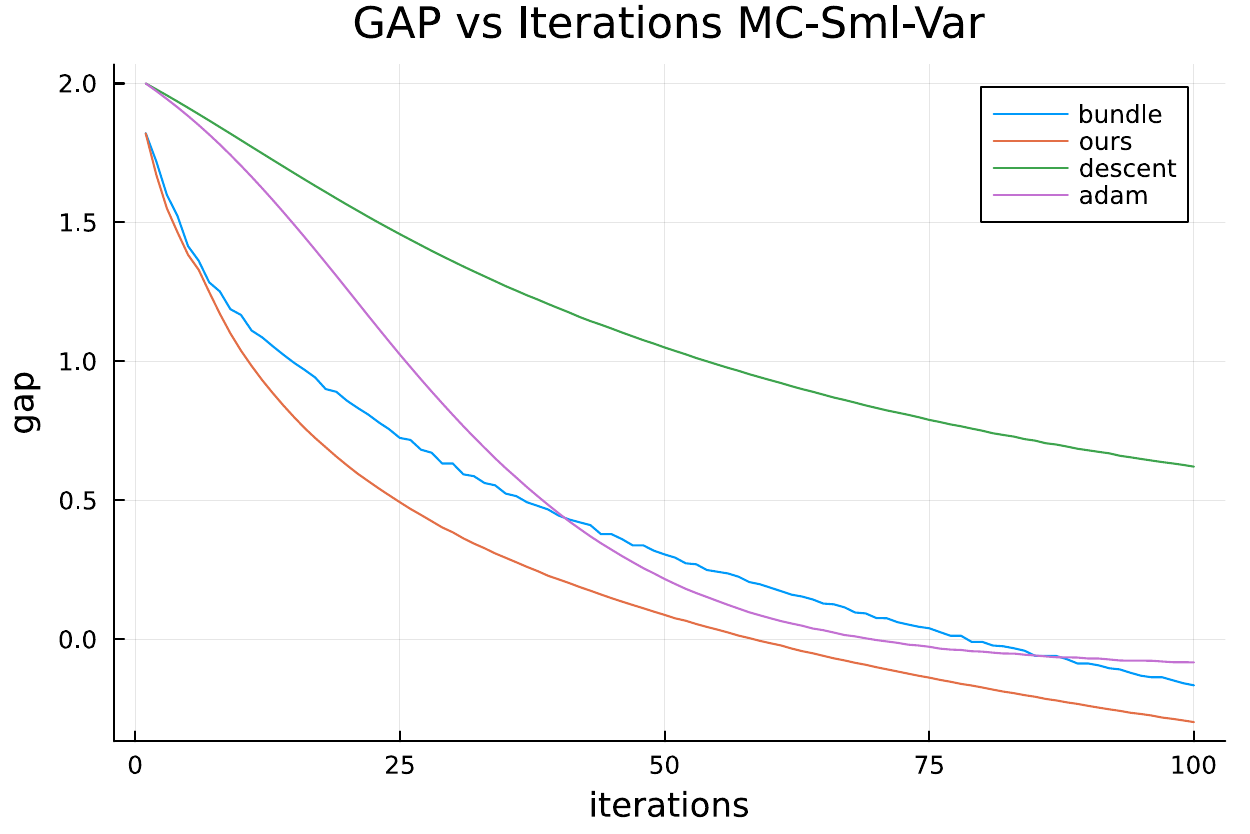}
   \includegraphics[width=0.45\textwidth]{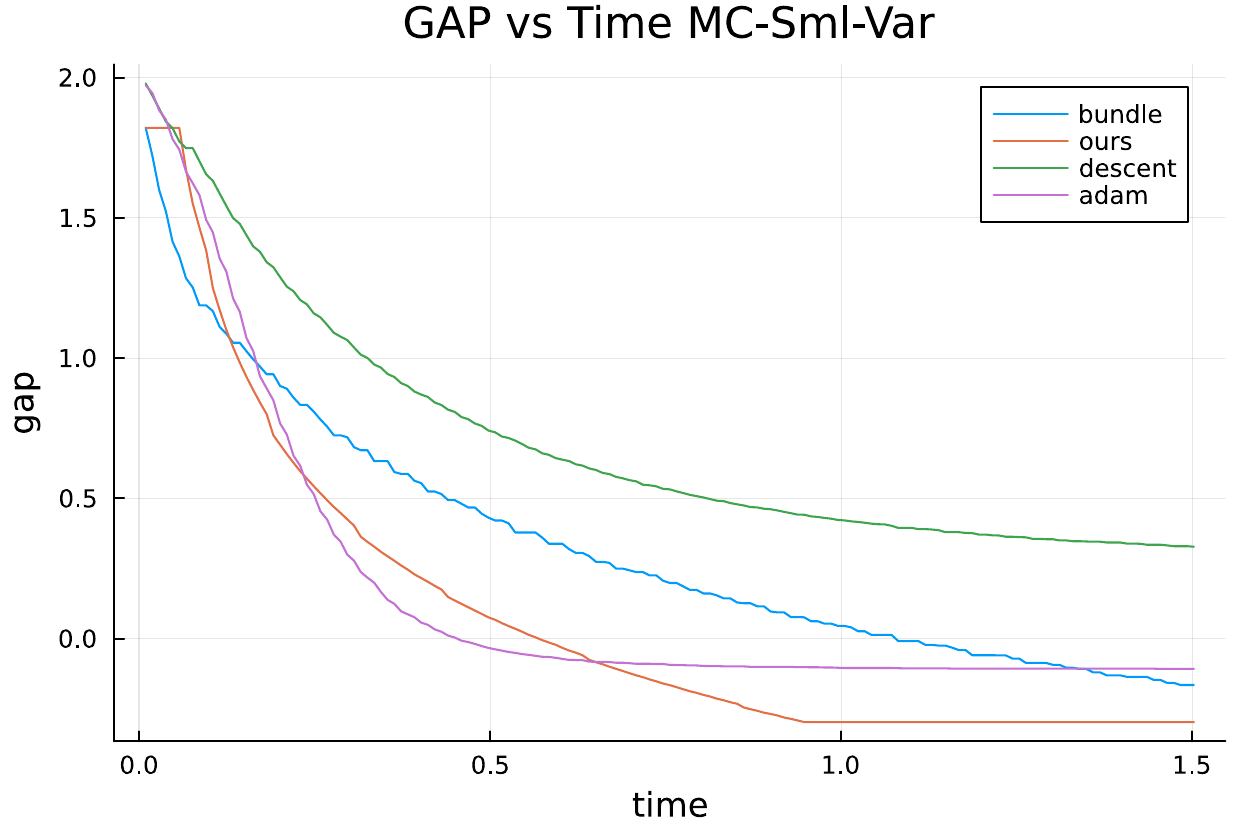}
    \caption{Comparison plots of the GAP (in logarithmic scale) during the execution of our trained \NN{} and different baselines in terms of iterations (in the left) and in terms of time (in the right) for \textsc{MC-Sml-Var} dataset.}
    \label{fig:execution_MC_Sml_Var}
\end{figure}

\begin{figure}
    \centering
   \includegraphics[width=0.45\textwidth]{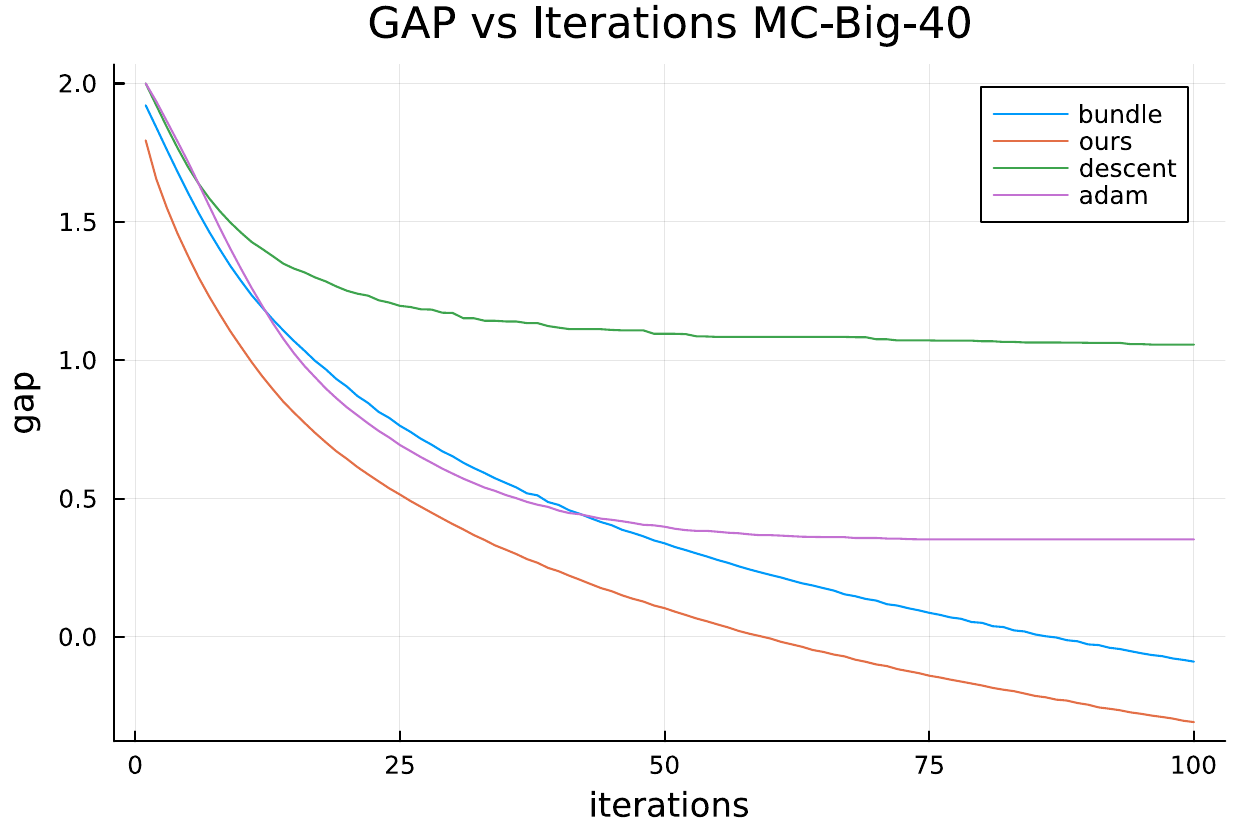}
   \includegraphics[width=0.45\textwidth]{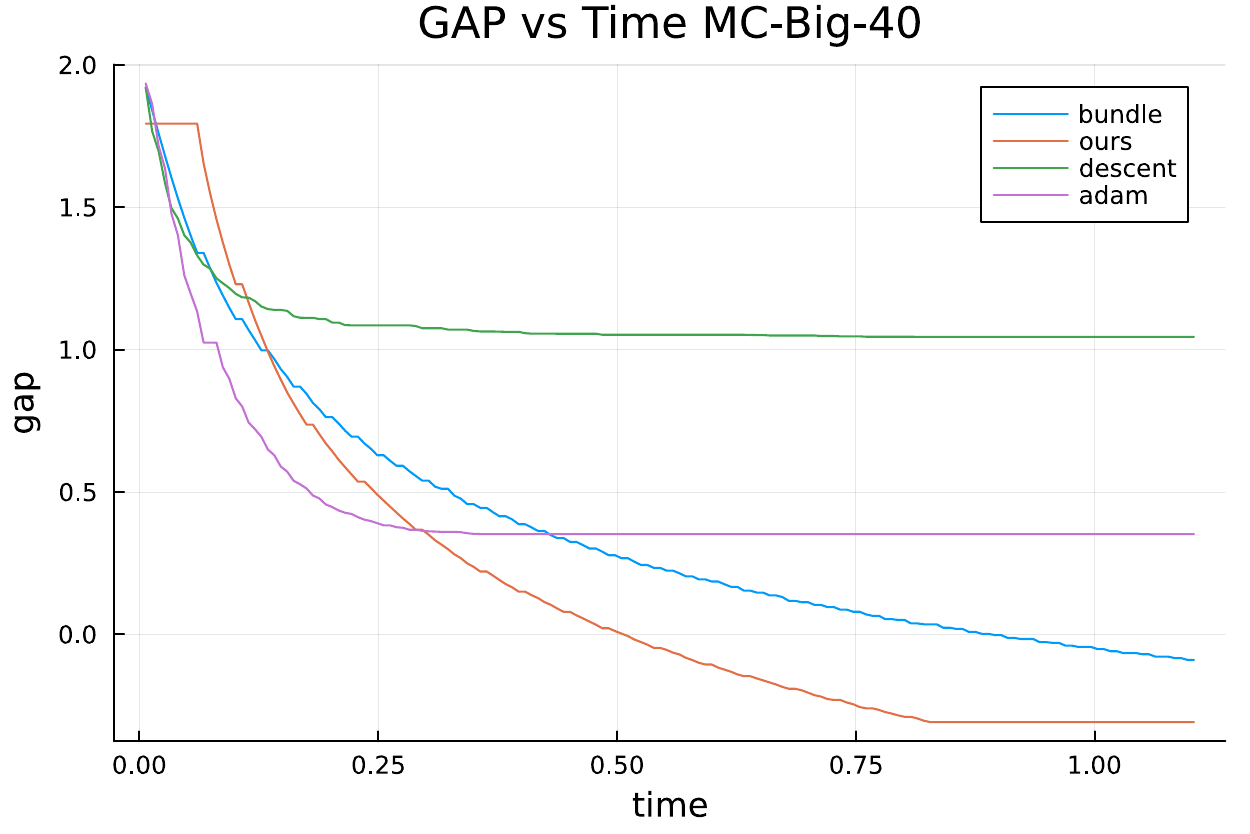}    
    \caption{Comparison plots of the GAP (in logarithmic scale) during the execution of our trained \NN{} and different baselines in terms of iterations (in the left) and in terms of time (in the right) for \textsc{MC-Big-40} dataset.}
    \label{fig:execution_MC_Big_40}
\end{figure}

\begin{figure}
    \centering
   \includegraphics[width=0.45\textwidth]{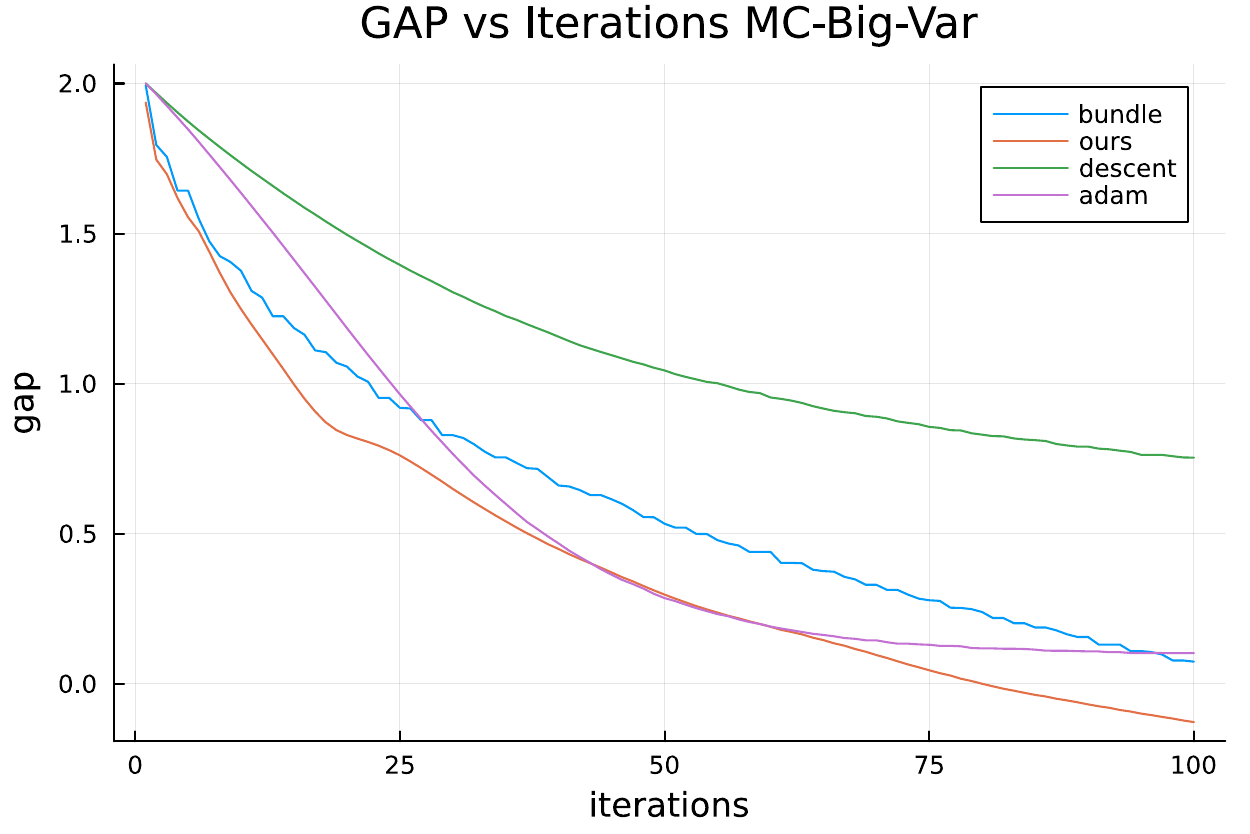}
   \includegraphics[width=0.45\textwidth]{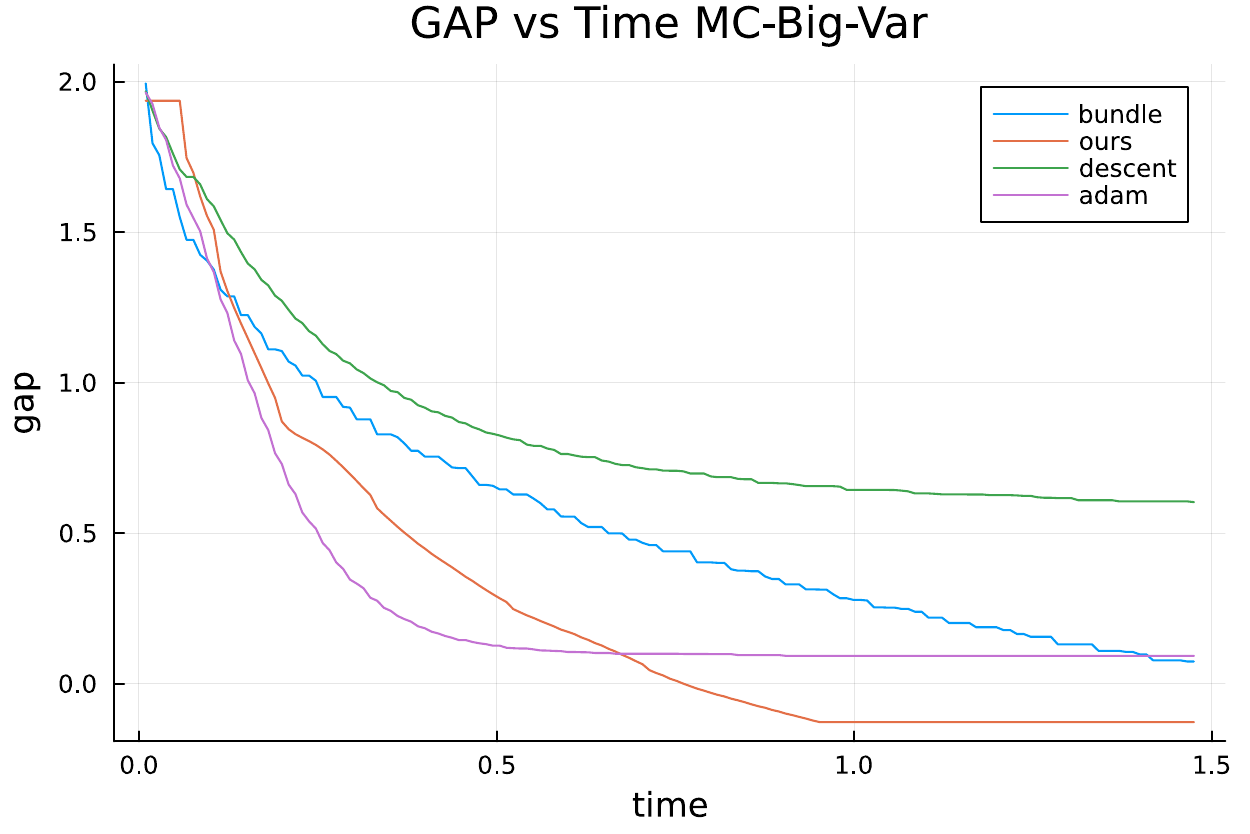}
    \caption{Comparison plots of the GAP (in logarithmic scale) during the execution of our trained \NN{} and different baselines in terms of iterations (in the left) and in terms of time (in the right) for \textsc{MC-Big-Var} dataset.}
    \label{fig:execution_MC_Big_Var}
\end{figure}

Anyway, by Figures~\ref{fig:execution_MC_Sml_40}-\ref{fig:execution_MC_Big_Var}, we can see that Bundle Network is slower than other approaches for the first iterations.
From these figures, we can see that initialization time can be improved.
These times can also depend on passages of CPU vectors to the GPU and vice versa.

From Table~\ref{tab:mainTable_bn} we can also see that the constant $\eta-$strategy provides almost every time the best performance, compared to the other \textit{deterministic} $\eta-$strategies.
This is a well-known characteristic of the \BM{}, as $100$ iterations is a \textit{relatively} small number of iterations for an \ABM{}, making $\eta$ tuning less performative than a constant strategy obtained after a grid-search.
This is because the \BM{} operates in distinct phases: initially, its goal is to identify promising directions for improvement, and later it focuses on collecting high-quality subgradients near the optimum to certify optimality.
This raises an interesting question about using multiple models to capture this Bundle behavior, which is beyond the scope of this work.

Nevertheless, our model provides high-quality predictions that generalize well to a larger number of iterations than those used during training (10).
In many cases, bundle network provides, in a given number of iterations, lower gaps than all the baselines, and even in cases where the provided predictions are not the best possibility, they lead to similar gaps to the optimal algorithm.

\begin{figure}
    \centering
    \centering
   \includegraphics[width=0.45\textwidth]{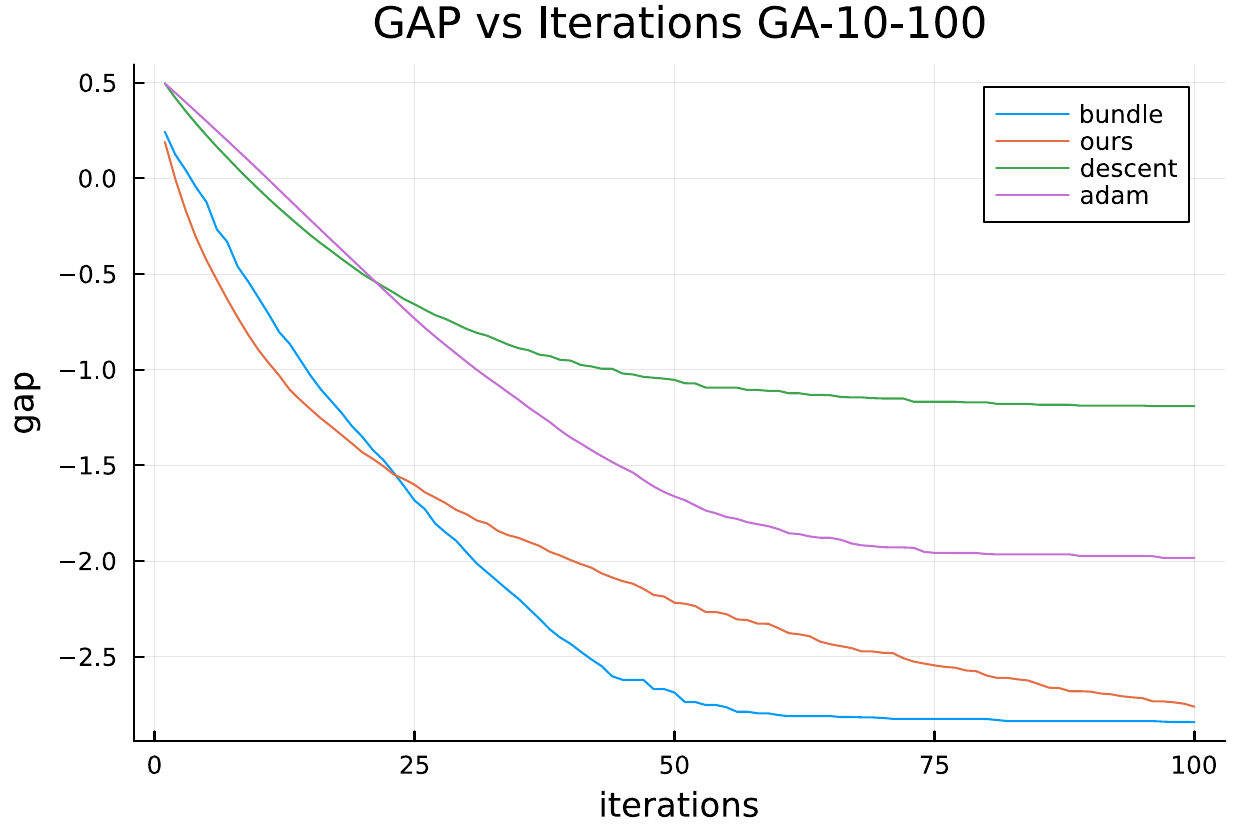}
   \includegraphics[width=0.45\textwidth]{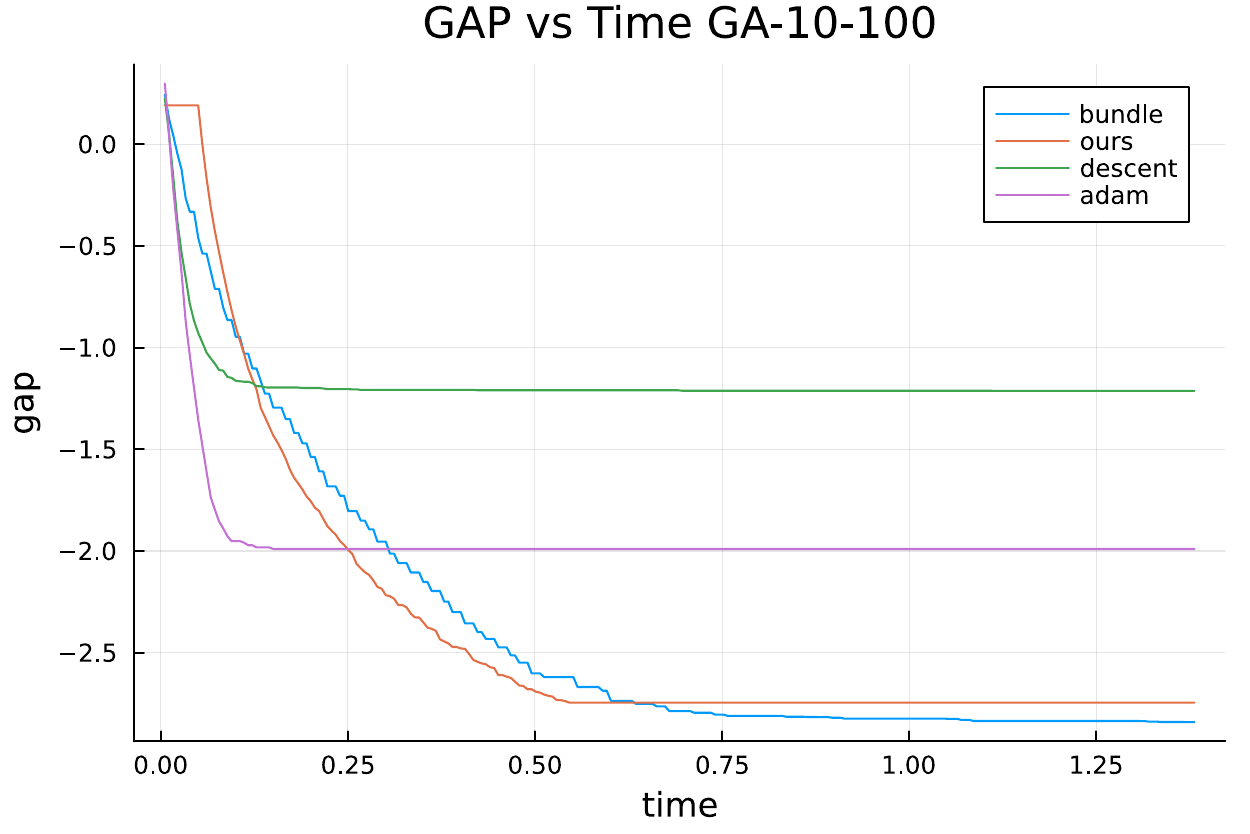}  
     \caption{Comparison plots of the GAP (in logarithmic scale) during the execution of our trained \NN{} and different baselines in terms of iterations (in the left) and in terms of time (in the right) for \textsc{GA-10-100} dataset.}
    \label{fig:execution_GA}
\end{figure}

\begin{figure}
    \centering
    \centering
   \includegraphics[width=0.45\textwidth]{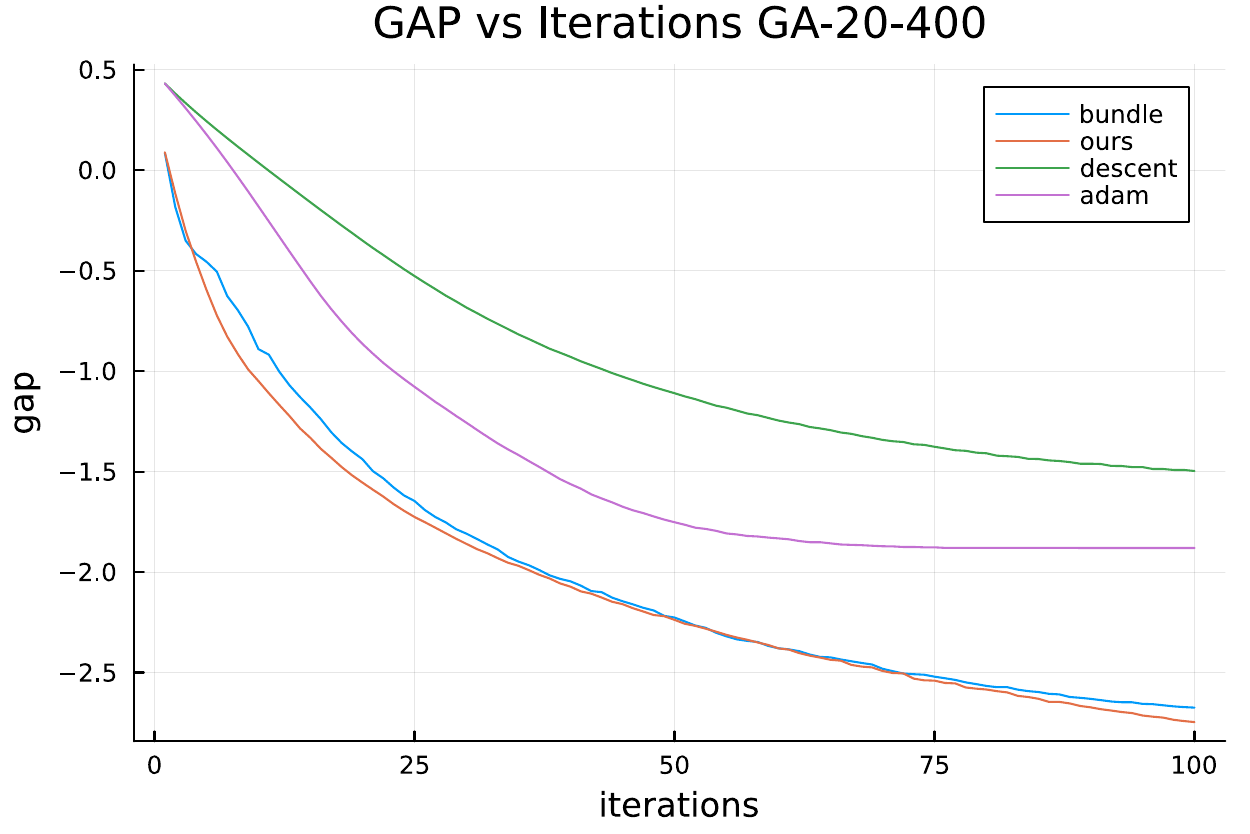}
   \includegraphics[width=0.45\textwidth]{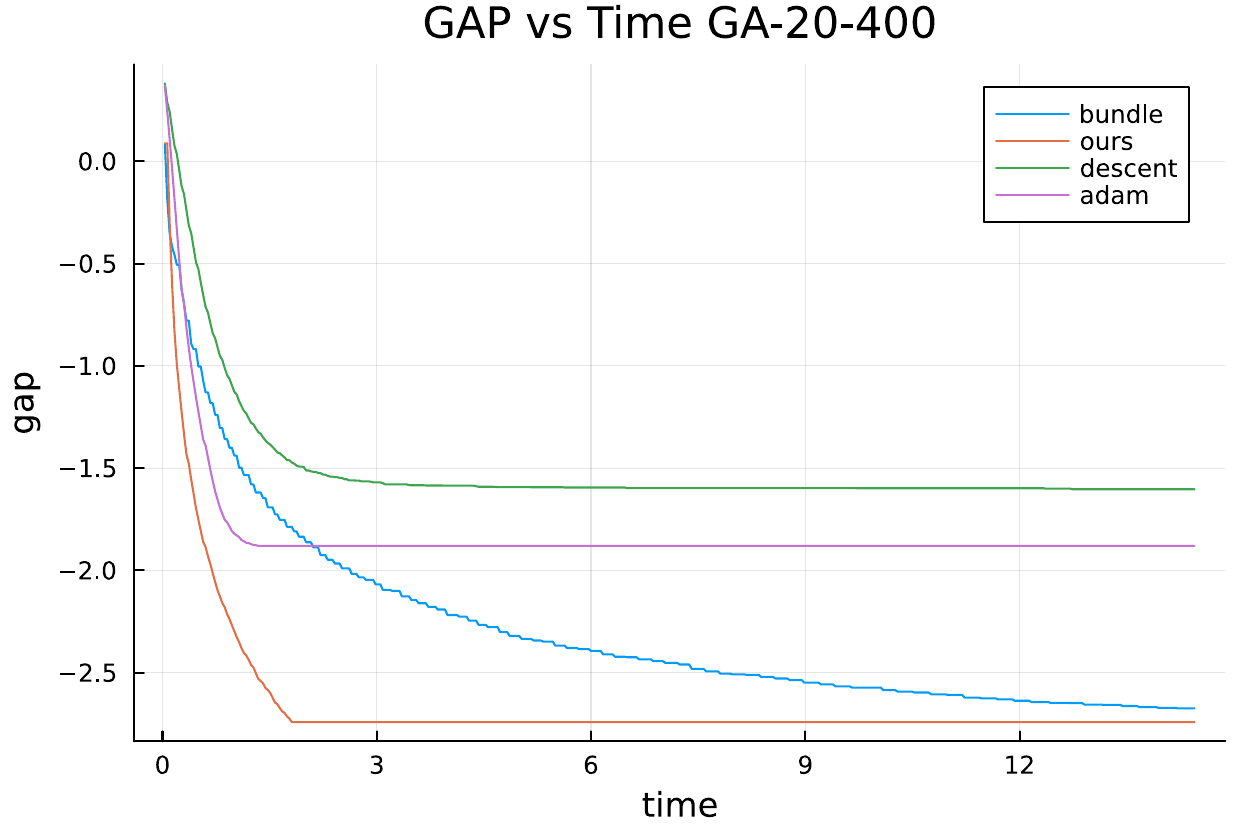}  
     \caption{Comparison plots of the GAP (in logarithmic scale) during the execution of our trained \NN{} and different baselines in terms of iterations (in the left) and in terms of time (in the right) for \textsc{GA-20-400} dataset.}
    \label{fig:execution_GA-20-400}
\end{figure}

\subsection{Sampling or Not Sampling}\label{sec:sampling}

\begin{table}[h]
\centering
\caption{Comparison of softmax and sparsemax, with or without using the sample strategy, or using the sample strategy only for t on MC-Sml40 }
    \label{tab:grid_search_MCsml40}
        \begin{tabular}{@{}lll rr c rr c rr c rr@{}} \toprule
Activation & \multicolumn{2}{c}{Sample } &  \multicolumn{2}{c}{10 iter.} & 
 \multicolumn{2}{c}{25 iter.} & 
& \multicolumn{2}{c}{50 iter.} & 
& \multicolumn{2}{c}{100 iter.} \\
\midrule
& t & q k
& GAP & time & 
& GAP & time & 
& GAP & time & 
& GAP & time \\
\midrule
softmax & & & \textbf{9.20} & 0.105 & & \textbf{1.75} & 0.185 & & \textbf{0.53} & 0.321 & & \textbf{0.17} & 0.592\\
softmax & x & & 9.78 & 0.104 & & 1.85 & 0.187 & & 0.56 & 0.321 & & \textbf{0.17} & 0.592\\
softmax & x & x & 10.41 & 0.105 & & 1.99 & 0.185 & & 0.60 & 0.318 & & 0.18 & 0.588\\
sparsemax & & & 58.3 & 0.125 & & 57.93 & 0.214 & & 57.47 & 0.361 & & 57.29 & 0.661\\
sparsemax & x & & 13.16 & 0.124 & & 2.61 & 0.212 & & 0.80 & 0.359 & & 0.24 & 0.655\\
sparsemax & x & x & 10.51 & 0.126 & & 2.07 & 0.214 & & 0.60 & 0.363 & & 0.18 & 0.662\\
\bottomrule
\end{tabular}
\end{table}

\begin{table}[h]
\centering
\caption{Comparison of softmax and sparsemax, with or without using the sample strategy, or using the sample strategy only for t on MC-Big-40 }
    \label{tab:grid_search_MCBig-40}
        \begin{tabular}{@{}lll rr c rr c rr c rr@{}} \toprule
Activation & \multicolumn{2}{c}{Sample } & \multicolumn{2}{c}{10 iter.} & 
& \multicolumn{2}{c}{25 iter.} & 
& \multicolumn{2}{c}{50 iter.} & 
& \multicolumn{2}{c}{100 iter.} \\
\midrule
& t & q k 
& GAP & time & 
& GAP & time & 
& GAP & time & 
& GAP & time \\
\midrule
softmax & & & \textbf{12.70} & 0.130 & & 3.41 & 0.243 & & 1.26 & 0.436 & & 0.47 & 0.832\\
softmax & x & & 12.85 & 0.131 & & \textbf{3.34} & 0.244 & & \textbf{1.22} & 0.436 & & \textbf{0.45} & 0.831\\
softmax & x & x & 12.88 & 0.131 & & 3.36 & 0.244 & & 1.26 & 0.436 & & 0.48 & 0.834\\
\midrule
sparsemax & & & 38.65 & 0.150 & & 24.00 & 0.270 & & 10.34 & 0.475 & & 3.39 & 0.898\\
sparsemax & x & & 39.40 & 0.150 & & 32.37 & 0.270 & & 17.76 & 0.474 & & 1.58 & 0.896\\
sparsemax & x & x & 13.27 & 0.150 & & 3.50 & 0.272 & & 1.40 & 0.479 & & 0.51 & 0.905\\
\bottomrule
\end{tabular}
\end{table}

\begin{table}[h]
\centering
\caption{Comparison of softmax and sparsemax, with or without using the sample strategy or using the sample strategy only for t on MC-Sml-Var }
    \label{tab:grid_search_MCsmlVar}
        \begin{tabular}{@{}lll rr c rr c rr c rr@{}} \toprule
Activation & \multicolumn{2}{c}{Sample } & & \multicolumn{2}{c}{10 iter.} 
& \multicolumn{2}{c}{25 iter.} & 
& \multicolumn{2}{c}{50 iter.} & 
& \multicolumn{2}{c}{100 iter.} \\
\midrule
& t & q k 
& GAP  & time & 
& GAP  & time & 
& GAP  & time & 
& GAP  & time  \\
\midrule
softmax &   &   &  \textbf{12.60} & 0.138 & & 3.28 & 0.269 & & \textbf{1.24} & 0.496 & & \textbf{0.50} & 0.952\\
softmax & x &   & 13.03 & 0.137 & & \textbf{3.27} & 0.268 & & 1.26 & 0.492 & & 0.51 & 0.948\\
softmax & x & x & 14.68 & 0.139 & & 3.94 & 0.269 & & 1.55 & 0.494 & & 0.62 & 0.957\\
\midrule
sparsemax &   &   & 30.23 & 0.157 & & 15.96 & 0.293 & & 12.45 & 0.527 & & 11.05 & 1.006\\
sparsemax & x &   & 56.83 & 0.160 & & 56.03 & 0.300 & & 56.00 & 0.538 & & 56.00 & 1.033\\
sparsemax & x & x & 14.02 & 0.158 & &  3.80 & 0.298 & &  1.56 & 0.531 & &  0.64 & 1.007\\
\bottomrule
\end{tabular}
\end{table}

\begin{table}[h]
\centering
\caption{Comparison of softmax and sparsemax, with or without using the sample strategy or using the sample strategy only for t on MC-Big-Var }
    \label{tab:grid_search_MCbigVar}
        \begin{tabular}{@{}lll rr c rr c rr c rr@{}} \toprule
Activation & \multicolumn{2}{c}{Sample } & \multicolumn{2}{c}{10 iter.} & 
& \multicolumn{2}{c}{25 iter.} & 
& \multicolumn{2}{c}{50 iter.} & 
& \multicolumn{2}{c}{100 iter.} \\
\midrule
& t & q k
& GAP & time & 
& GAP & time & 
& GAP & time & 
& GAP & time \\
\midrule
softmax &   &   & 20.10 & 0.139 & & 5.11 & 0.271 & & 2.01 & 0.495 & & 0.70 & 0.953\\
softmax & x &  & 16.79 & 0.141 & & 4.78 & 0.273 & & 1.68 & 0.499 & & 0.63 & 0.956\\
softmax & x & x & \textbf{16.66} & 0.139 & & 4.63 & 0.270 & & 1.80 & 0.495 & & 0.67 & 0.955\\
\midrule
sparsemax & & & 42.84 & 0.160 & & 32.39 & 0.301 & & 29.68 & 0.541 & & 23.25 & 1.031\\
sparsemax & x & & 30.20 & 0.162 & & 12.33 & 0.305 & & 6.25 & 0.542 & & 3.78 & 1.033\\
sparsemax & x & x & 18.94 & 0.162 & & \textbf{4.24} & 0.305 & & \textbf{1.47} & 0.548 & & \textbf{0.61} & 1.045\\
\bottomrule
\end{tabular}
\end{table}

\begin{table}[h]
\centering
\caption{Comparison of softmax and sparsemax, with or without using the sample strategy, or using the sample strategy only for t on GA-10-100 }
    \label{tab:grid_search_GA}
        \begin{tabular}{@{}lll rr c rr c rr c rr@{}} \toprule
Activation & \multicolumn{2}{c}{Sample } & \multicolumn{2}{c}{10 iter.} & 
& \multicolumn{2}{c}{25 iter.} & 
& \multicolumn{2}{c}{50 iter.} & 
& \multicolumn{2}{c}{100 iter.} \\
\midrule
& t & q k& 
GAP & time & 
& GAP & time & 
& GAP & time & 
& GAP & time \\
\midrule
softmax & & & \textbf{0.1484} & 0.104 & & \textbf{0.0228} & 0.177 & & \textbf{0.0047} & 0.304 & & \textbf{0.0009} & 0.551\\
softmax & x & & 0.1544 & 0.102 & & 0.0229 & 0.176 & & 0.0049 & 0.303 & & \textbf{0.0009} & 0.548\\
softmax & x & x & 0.1954 & 0.099 & & 0.0370 & 0.176 & & 0.0081 & 0.302 & & 0.0013 & 0.551\\
\midrule
sparsemax & & & 0.2335 & 0.118 & & 0.0600 & 0.201 & & 0.0336 & 0.339 & & 0.0288 & 0.615\\
sparsemax & x & & 1.0482 & 0.120 & & 0.9783 & 0.202 & & 0.9782 & 0.342 & & 0.9782 & 0.611\\
sparsemax & x & x & 0.9172 & 0.122 & & 0.7001 & 0.201 & & 0.6827 & 0.340 & & 0.6818 & 0.609\\
\bottomrule
\end{tabular}
\end{table}

\begin{table}[h]
\centering
\caption{Comparison of softmax and sparsemax, with or without using the sample strategy, or using the sample strategy only for t on GA-20-400 }
    \label{tab:grid_search_GA-20-400}
        \begin{tabular}{@{}lll rr c rr c rr c rr@{}} \toprule
Activation & \multicolumn{2}{c}{Sample } & \multicolumn{2}{c}{10 iter.} & 
& \multicolumn{2}{c}{25 iter.} & 
& \multicolumn{2}{c}{50 iter.} & 
& \multicolumn{2}{c}{100 iter.} \\
\midrule
& t & q k
& GAP & time & 
& GAP & time & 
& GAP & time & 
& GAP & time \\
\midrule
softmax & & & \textbf{0.1018} & 0.228 & & \textbf{0.0190} & 0.498 & & \textbf{0.0052} & 0.943 & & \textbf{0.0014} & 1.837\\
softmax & x & & 0.1065 & 0.218 & & 0.0216 & 0.480 & & 0.0060 & 0.904 & & 0.0014 & 1.756\\
softmax & x & x & 0.1273 & 0.204 & & 0.0295 & 0.445 & & 0.0087 & 0.836 & & 0.0018 & 1.621\\
\midrule
sparsemax & & & 0.1463 & 0.235 & & 0.0626 & 0.500 & & 0.0601 & 0.935 & & 0.0599 & 1.804\\
sparsemax & x & & 0.7860 & 0.219 & & 0.7696 & 0.463 & & 0.767 & 0.860 & & 0.7669 & 1.656\\
sparsemax & x & x & 0.4002 & 0.218 & & 0.1272 & 0.461 & & 0.0681 & 0.860 & & 0.0627 & 1.651\\
\bottomrule
\end{tabular}
\end{table}

In this section, we compare the model trained with and without the sampling strategy, for the hidden representation discussed earlier, as well as considering two different possibilities for $\psi$ that is, the softmax~\citep[pp. 180-184]{goodfellow2016deep} or the sparsemax~\citep{MartinsA16}.
In all the cases the models are trained unrolling 10 iterations of the resolutive approach and for $25$ training epochs using Adam as optimizer, with learning rate $10^{-5}$ and a Clip Norm (to 5).
.
In the case in which we do not use the sampling mechanism, the hidden representations for the queries, the key, and the step size are not predicted by sampling them through a Gaussian distribution. Indeed, we learn the mean and the variance.
Instead, it is predicted directly by the encoder.
In the test phase, the hidden representation is not sampled for either approach. Instead, we take the mean directly when the model predicts both the mean and the variance.

From Tables \ref{tab:grid_search_MCsml40}-\ref{tab:grid_search_GA-20-400}, we can see that sampling is a strategy leading to better performance when using the sparsemax, but it is inefficient for the softmax. This can be because the sparsemax function has a zero gradient in a certain region of the space that can be escaped using sampling strategies.

We choose the softmax with no sampling strategy as it shows better performances that are stable across both problems.
Meanwhile, the sparsemax seems to be an interesting choice for harder datasets, as shown in Table \ref{tab:grid_search_MCbigVar}.

\subsection{Testing on another dataset}\label{sec:cross_testing}

From Tables~\ref{tab:gen_s40}-\ref{tab:gen_ga2}, we can see the performance of a model trained on one dataset and tested on another.
Training the model on another dataset rarely provides better gaps than training it on the same dataset. It happens in Table \ref{tab:gen_s40} and for 10 iterations in Table \ref{tab:gen_ga1}.
Considering the same problem, we can obtain similar behaviors in several cases, in particular while using more iterations.
In some cases, the model trained on the GA dataset achieves interesting performance on the small datasets of MC, as shown in Tables \ref{tab:gen_s40} and \ref{tab:gen_sVar}.
Anyway, no update learned for the MC dataset seems to work sufficiently well for the GA datasets.
Tables~\ref{tab:gen_s40}-\ref{tab:gen_ga2} also show that, in several cases, it is possible to attain lower gaps for iteration from the performances of a model trained on another dataset, of the same problem, than the baseline, providing a lower gap obtained with grid search.

\begin{table}[h!]
\caption{Gap of the model tested on MC-Sml-40 trained on different datasets. The associated training dataset is displayed as the first element of each row of the table. The last row shows the performances of the best baseline obtained with grid search.}
\label{tab:gen_s40}
\begin{tabular}{@{}l c rr c rr c rr c rr@{}}
\toprule
\multirow{2}{*}{Training Dataset} & & \multicolumn{2}{c}{10 iter.} & 
& \multicolumn{2}{c}{25 iter.} & 
& \multicolumn{2}{c}{50 iter.} & 
& \multicolumn{2}{c}{100 iter.} \\
\cmidrule{3-4} \cmidrule{6-7} \cmidrule{9-10} \cmidrule{12-13}
& & GAP & time & 
& GAP & time & 
& GAP & time & 
& GAP & time \\
\midrule
\textbf{MC-Sml-40 }& & 9.20 & 0.105 & & 1.75 & 0.185 & & 0.53 & 0.321 & & \textbf{0.17} & 0.592\\
MC-Sml-Var & & \textbf{8.02} & 0.104 & & \textbf{1.69} & 0.184 & & \textbf{0.51} & 0.319 & & \textbf{0.17} & 0.593\\
MC-Big-40  & & 18.02 & 0.116 & & 4.14 & 0.206 & & 1.49 & 0.358 & & 0.47 & 0.665\\
MC-Big-Var & & 10.77 & 0.106 & & 2.26 & 0.187 & & 0.70 & 0.327 & & 0.20 & 0.603\\
\midrule
GA-10-100 & & 61.98 & 0.104 & & 20.79 & 0.183 & &  4.07 & 0.319 & & 0.95 & 0.589\\
GA-20-400 & & 86.33 & 0.107 & & 68.16 & 0.190 & & 14.20 & 0.327 & & 3.52 & 0.607\\
\midrule
Best Baseline & & 11.96 & 0.048 & & 4.16 & 0.128 & & 1.30 & 0.293 & & 0.31 & 0.816\\
\bottomrule
     \end{tabular}
\end{table}

\begin{table}[h!]
\caption{Gap of the model tested on MC-Sml-Var trained on different datasets. The associated training dataset is displayed as first element of each row of the table. The last row shows the performances of the best baseline obtained with grid-search.}
\label{tab:gen_sVar}
\begin{tabular}{@{}l c rr c rr c rr c rr@{}}
\toprule
\multirow{2}{*}{Training Dataset} & & \multicolumn{2}{c}{10 iter.} & 
& \multicolumn{2}{c}{25 iter.} & 
& \multicolumn{2}{c}{50 iter.} & 
& \multicolumn{2}{c}{100 iter.} \\
\cmidrule{3-4} \cmidrule{6-7} \cmidrule{9-10} \cmidrule{12-13}
& & GAP & time & 
& GAP & time & 
& GAP & time & 
& GAP & time \\
\midrule
MC-Sml-40 & & 19.87 & 0.142 & & 5.72 & 0.276 & & 2.32 & 0.501 & & 1.00 & 0.965\\
\textbf{MC-Sml-Var} & & \textbf{12.60} & 0.138 & & \textbf{3.28} & 0.269 & & \textbf{1.24} & 0.496 & & \textbf{0.50} & 0.952\\
MC-Big-40 & & 33.83 & 0.149 & & 14.03 & 0.291 & & 7.54 & 0.536 & & 3.90 & 1.033\\
MC-Big-Var & & 16.21 & 0.141 & & 3.83 & 0.276 & & 1.39 & 0.505 & & 0.52 & 0.976\\
\midrule
GA-10-100 & & 25.76 & 0.142 & & 7.60  & 0.275 & &  2.08 & 0.502 & &  0.67 & 0.969\\
GA-20-400 & & 66.15 & 0.140 & & 62.06 & 0.275 & & 38.19 & 0.507 & & 12.86 & 0.974\\
\midrule
Best Baseline & & {12.79} & 0.108 &  & 4.33 & 0.29 & & 1.66  &  0.635 & & 0.58 & 1.503\\
\bottomrule
\end{tabular}
\end{table}

\begin{table}[h!]
\caption{Gap of the model tested on MC-Big-40 trained on different datasets. The associated training dataset is displayed as first element of each row of the table. The last row shows the performances of the best baseline obtained with grid-search.}
\begin{tabular}{@{}l c rr c rr c rr c rr@{}}
\toprule
\multirow{2}{*}{Training Dataset} & & \multicolumn{2}{c}{10 iter.} & 
& \multicolumn{2}{c}{25 iter.} & 
& \multicolumn{2}{c}{50 iter.} & 
& \multicolumn{2}{c}{100 iter.} \\
\cmidrule{3-4} \cmidrule{6-7} \cmidrule{9-10} \cmidrule{12-13}
& & GAP & time & 
& GAP & time & 
& GAP & time & 
& GAP & time \\
\midrule
MC-Sml-40 & & 35.09 & 0.133 & & 5.78 & 0.248 & & 1.63 & 0.445 & & 0.62 & 0.85\\
MC-Sml-Var & & 56.68 & 0.132 & & 25.62 & 0.246 & & 8.88 & 0.441 & & 3.06 & 0.842\\
\textbf{MC-Big-40} & & \textbf{12.70} & 0.130 & & \textbf{3.41} & 0.243 & & \textbf{1.26} & 0.436 & & \textbf{0.47} & 0.832\\
MC-Big-Var & & 37.57 & 0.134 & & 10.74 & 0.251 & & 4.89 & 0.449 & & 1.49 & 0.856\\
\midrule
GA-10-100 & & 95.70 & 0.123 & & 95.70 & 0.230 & & 78.63 & 0.408 & & 27.53 & 0.779\\
GA-20-400 & & 100.0 & 0.134 & & 100.0 & 0.252 & & 99.53 & 0.447 & & 74.54 & 0.854\\
\midrule
Best Baseline & & 19.34 & 0.076 & & 5.64 & 0.203 & & 2.08 &  0.447 & & 0.78 & 1.110\\
\bottomrule
     \end{tabular}
\end{table}

\begin{table}[h!]
\caption{Dataset tested on MC-Big-Var trained on different datasets. The associated training dataset is displayed as first element of each row of the table. The last row shows the performances of the best baseline obtained with grid-search.}
\begin{tabular}{@{}l c rr c rr c rr c rr@{}}
\toprule
\multirow{2}{*}{Training Dataset} & & \multicolumn{2}{c}{10 iter.} & 
& \multicolumn{2}{c}{25 iter.} & 
& \multicolumn{2}{c}{50 iter.} & 
& \multicolumn{2}{c}{100 iter.} \\
\cmidrule{3-4} \cmidrule{6-7} \cmidrule{9-10} \cmidrule{12-13}
& & GAP & time & 
& GAP & time & 
& GAP & time & 
& GAP & time \\
\midrule
MC-Sml-40 & & 21.50  & 0.142 & & 5.48  & 0.280 & & 2.07 & 0.511 & & 0.85 & 0.984\\
MC-Sml-Var & & 21.03 & 0.221 & & 7.64  & 0.419 & & 2.67 & 0.755 & & 0.95 & 1.430\\
MC-Big-40 & & 28.39  & 0.150 & & 11.15 & 0.294 & & 5.87 & 0.540 & & 2.94 & 1.037\\
\textbf{MC-Big-Var} & & {20.10} & 0.139 & & \textbf{5.11} & 0.271 & & \textbf{2.01} & 0.495 & & \textbf{0.70} & 0.953\\
\midrule
GA-10-100 & & 41.66 & 0.144 & & 28.58 & 0.280 & & 20.72 & 0.511 & & 8.30 & 0.979\\
GA-20-400 & & 72.89 & 0.144 & & 69.99 & 0.283 & & 53.10 & 0.519 & & 27.75 & 0.996\\
\midrule
Best Baseline & & \textbf{17.98} & 0.112 & & 6.34 & 0.289 &  & 2.66 & 0.627 & & 0.93 & 1.477\\
\bottomrule     \end{tabular}
\end{table}

\begin{table}[h!]
\caption{Gap of the model tested on GA-10-100 trained on different datasets. The associated training dataset is displayed as first element of each row of the table. The last row shows the performances of the best baseline obtained with grid-search.}\label{tab:gen_ga1}
\begin{tabular}{@{}l c rr c rr c rr c rr@{}}
\toprule
\multirow{2}{*}{Training Dataset} & & \multicolumn{2}{c}{10 iter.} & 
& \multicolumn{2}{c}{25 iter.} & 
& \multicolumn{2}{c}{50 iter.} & 
& \multicolumn{2}{c}{100 iter.} \\
\cmidrule{3-4} \cmidrule{6-7} \cmidrule{9-10} \cmidrule{12-13}
& & GAP & time & 
& GAP & time & 
& GAP & time & 
& GAP & time \\
\midrule
MC-Sml-40  & & 2.5903 & 0.101 & & 1.8558 & 0.179 & & 0.9291 & 0.304 & & 0.4860 & 0.555\\
MC-Sml-Var & & 2.5165 & 0.102 & & 1.4663 & 0.178 & & 1.0522 & 0.307 & & 1.0512 & 0.562\\
MC-Big-40  & & 2.7488 & 0.103 & & 2.1948 & 0.177 & & 1.3515 & 0.302 & & 0.7463 & 0.550\\
MC-Big-Var & & 2.4958 & 0.102 & & 1.3267 & 0.180 & & 0.9705 & 0.309 & & 0.9683 & 0.564\\
\midrule
\textbf{GA-10-100} & & 0.1484 & 0.104 & & {0.0228} & 0.177 & & {0.0048} & 0.304 & & {0.0009} & 0.551\\
GA-20-400 & & \textbf{0.1332} & 0.105 & & 0.0233 & 0.182 & & 0.0056 & 0.309 & & 0.0016 & 0.561\\
\midrule
Best Baseline & & {0.1893} & 0.071 & & \textbf{0.0156} & 0.211 & & \textbf{0.0014} & 0.601 & & \textbf{0.0006} & 1.596\\
\bottomrule
     \end{tabular}
\end{table}

\begin{table}[h!]
\caption{Gap of the model tested on GA-20-400 trained on different datasets. The associated training dataset is displayed as first element of each row of the table. The last row shows the performances of the best baseline obtained with grid-search.}\label{tab:gen_ga2}
\begin{tabular}{@{}l c rr c rr c rr c rr@{}}
\toprule
\multirow{2}{*}{Training Dataset} & & \multicolumn{2}{c}{10 iter.} & 
& \multicolumn{2}{c}{25 iter.} & 
& \multicolumn{2}{c}{50 iter.} & 
& \multicolumn{2}{c}{100 iter.} \\
\cmidrule{3-4} \cmidrule{6-7} \cmidrule{9-10} \cmidrule{12-13}
& & GAP & time & 
& GAP & time & 
& GAP & time & 
& GAP & time \\
\midrule
MC-Sml-40  & & 2.3695 & 0.221 & & 1.9535 & 0.474 & & 1.1906 & 0.894 & & 0.5052 & 1.732\\
MC-Sml-Var & & 2.3357 & 0.231 & & 1.7168 & 0.507 & & 0.8963 & 0.959 & & 0.8654 & 1.868\\
MC-Big-40  & & 2.4715 & 0.218 & & 2.1550 & 0.481 & & 1.5618 & 0.909 & & 0.8056 & 1.759\\
MC-Big-Var & & 2.3197 & 0.218 & & 1.6304 & 0.482 & & 0.8092 & 0.910 & & 0.7856 & 1.765\\
\midrule
GA-10-100 & & 0.1811 & 0.234 & & 0.0296 & 0.511 & & 0.0077 & 0.962 & & 0.0015 & 1.866\\
\textbf{GA-20-400} & & \textbf{0.1018} & 0.228 & & \textbf{0.0190} & 0.498 & & 0.0052 & 0.943 & & \textbf{0.0014} & 1.837\\
\midrule
Best Baseline & & 0.1129 & 0.433 & & 0.0193 & 1.459 & & \textbf{0.0048} & 4.256 & & 0.0020 & 14.459\\
\bottomrule
     \end{tabular}
\end{table}

\section{Conclusion}

In this work, we present a machine learning model based on unrolling and amortization strongly inspired by the bundle method. The model bypasses the need to solve the \DMP{} at each iteration, leading to a smooth optimization process that can be differentiated using Automatic Differentiation techniques. 
We test it for the resolution of the Lagrangian dual obtained from the two different problems: the Multi-Commodity Network Design and the Generalized Assignment.
It achieves a lower gap compared to other classical approaches for a fixed number of iterations.
Moreover, it exhibits strong generalization capabilities, maintaining robust performance even beyond the training horizon and sometimes across different datasets.

This work opens several promising directions for future research.
First, the framework can be directly applied, as it is, to other problems arising from Lagrangian relaxation. A central challenge in standard bundle methods lies in automatically determining which components to retain or discard in order to control bundle size and improve computational efficiency. While straightforward modifications of the framework are possible, this also raises questions concerning feature extraction and the design of an automated feature selection phase.

Given the growing use of bundle methods in machine learning, it is also interesting to assess the performance of the proposed approach in meta-learning tasks. This requires applying the method to non-convex and non-smooth optimization problems and this extension presents non-trivial challenges. In this setting, removing components from the bundle is not only essential for reducing computational cost, but also for avoiding convergence to poor local minima.

\backmatter

\bibliography{biblio}

\begin{appendices}

\section{Proximality Term $\eta$-Strategy}\label{app:eta_strat}

The proper selection of the proximal term $\eta$ of the bundle method is both non-trivial and essential.
This choice depends on the current stabilization point and the information accumulated in the bundle up to the current iteration. 
In particular, when all gradients at the optimum are available, solving a \CP{} model is preferable, as it yields the optimal solution to the problem.
Hence, in this case, a higher value of $\eta$ may be preferred.
While at the beginning of the resolution we possibly dispose of very local information, acting carefully is preferred, and so a smaller value of $\eta$ may be preferred.

In \citep{frangioni2002generalized}, theoretical conditions are established for the $\eta$-strategy to ensure convergence.
Constant strategies can be effective if appropriately fine-tuned.
Many state-of-the-art strategies are self-adjusting rules for choosing $\eta_t$.

The main state-of-the-art $\eta$-strategies are composed of three \textit{levels}, as the ones developed in \citep{frangioni1997dual} and also presented in \citep{crainic_bundle-based_2001}.
These levels includes: \textit{Long-term}, \textit{Middle-term}, and \textit{Short-term} $\eta$-\textit{strategies}.

The \textit{Long-term $\eta$-strategies} aim to capture the stage of the resolution process. In the \BM{}, this is important because it represents a balance between the \SGM{} and the \CPM{}.
In particular, we want to prioritize \SGM{} iterations when the solution is still far from optimal, and shift towards \CP{} in the final iterations.
This approach helps collect the subgradient information necessary to establish the optimality of the solution.

We define the maximum improvement estimation of the \DMP{} as
$$
v^*_t=\eta_t ||\sum_{i\in\beta_t}\theta_{i}^{(t)}\vg_i||^2 + \sum_{i\in\beta_t}\alpha_i\theta_{i}^{(t)}
$$
where $\vtheta^{(i)}$ is the solution of the \DMP{} at iteration $t$.

We also define the expected maximum improvement obtained using a variant of the \DMP{} in which $\eta^*$ acts as a regularization parameter:
$$
\epsilon^*_t = \sum_{i\in\beta_t}\alpha_i\theta_{i}^{(t)} + \eta^* ||\sum_{i\in\beta_t}\theta_{i}^{(t)}\vg_i||^2
$$
where $\eta^*$ is a hyperparameter.
The hyperparameter $\eta^*$ is generally chosen to be greater than all the possible $\eta_t$, thereby emulating the expected behavior by following \textit{more} the decision of the model.
It is different from explicitly considering this parameter in the \DMP{}, as this would alter its solution and consequently the new trial direction, as previously explained.
Additionally, we define a hyperparameter $\Tilde{m} \in [0, 1)$, which is employed in the $\eta$-strategies and is usually fixed to $0.01$.

Some examples of Long-term $\eta$-strategies are:
\begin{itemize}
    \item no Long-term $\eta$-strategy;
    \item  the \textit{soft} Long-term $\eta$-strategy. After a null step, decreases in $\eta_t$ are inhibited whenever
                $ v^*_t < \Tilde{m} \epsilon^*_t$.
    In other words, decrements are prevented if $\eta_t$ is already sufficiently small, such that the maximum improvement is estimated to be less than $\Tilde{m}$ times the improvement suggested by a Cutting-Planes model.
            
    \item the \textit{hard} Long-term $\eta$-strategy. Increases of $\eta_t$ are allowed, even after null steps, whenever
                $ v^*_t < \Tilde{m}  \epsilon^*_t$.
    
    \item the \textit{balancing} Long-term $\eta$-strategy. This strategy aims to keep the two
              terms $\eta^*||\frac{1}{2}\sum_{i\in\beta_t}\theta_{i}^{(t)}\vg_i||^2_2$ {roughly the same size}.
             
              Increases of $\eta_t$, after a serious step, 
             are inhibited if 
             $$\frac{\eta^*}{2}||\sum_{i\in\beta_t}\theta_{i}^{(t)}\vg_i||^2_2 \leq \Tilde{m}  \sum_{i=1}^{|\beta_t|}\alpha_i\theta_i^{(t)}.
             $$
             In this case, increasing $\eta_t$ causes a reduction of $\frac{\eta^*}{2}||\sum_{i\in\beta_t}\theta_{i}^{(t)}\vg_i||^2_2$ 	         which is already small.
           
             On the other hand, decreases in $\eta_t$  are inhibited if 
             $$
             \Tilde{m}  \frac{\eta^*}{2}||\sum_{i\in\beta_t}\theta_{i}^{(t)}\vg_i||^2_2 \geq \sum_{i=1}^{|\beta_t|}\alpha_i\theta_i^{(t)}.
             $$
             The rational been that decreasing $\eta_t$ causes an
             increase in $||\sum_{i\in\beta_t}\theta_{i}^{(t)}\vg_i||^2_2$ which is already big.
\end{itemize}

The \textit{Middle-term $\eta$-strategies} are multi-step strategies that consider the outcomes from the \textit{last few iterations}.
These strategies include the following conditions:

\begin{itemize}
\item A minimum number of consecutive successful steps (SS) with the same $\eta_t$ must be performed before $\eta_t$ is allowed to increase. 
\item A minimum number of consecutive null steps (NS) with the same $\eta_t$ must be performed before $\eta_t$ is allowed to decrease. 
\end{itemize}

At the bottom, heuristic \textit{Short-term $\eta$-strategies} only rely on information from the current iteration to propose an update for $\eta_t$. Specifically, these heuristics are only used to assist in selecting the current value, but both the top and the middle layers would agree on the update. 
 However, sometimes it can be beneficial to allow the Short-term heuristics to do \textit{small adjustments} in $\eta_t$, regardless of the decisions made by the other levels.

We consider in our implementation these \textit{Short-term $\eta$-strategies} together :
\begin{itemize}
    \item  Increase $\eta_{t+1}$ as $\eta_{t} \eta_{incr}$, with $\eta_{incr}>1$. 
\item Decrease $\eta_{t+1}$ as $\eta_{t-1}  \eta_{decr}$, with $0<\eta_{decr}<1$. 
\item When increased, $\eta_{t+1}$ should not exceed $\eta_M$, with $\eta_M>0$.
\item When decreased, $\eta_{t+1}$ should not go below $\eta_m$, with $0<\eta_m < \eta_M$. 
\end{itemize}

\section{Multi Commodity Capacitated Network Design Problem \label{app:MC}}

A MC instance is given by a directed simple graph $D = (N, A)$, a set of commodities $K$, an arc-capacity vector $c$, and two cost vectors $r$ and $f$.
Each commodity $k\in K$ corresponds to a triplet $(o^k, d^k, q^k)$ where $o^k\in N$ and $d^k\in N$ are the nodes corresponding to the origin and the destination of commodity $k$, and $q^k \in \mathbb{N}^*$ is its volume. For each arc, $(i,j) \in A$, $c_{ij}>0$ corresponds to the maximum amount of flow that can be routed through $(i,j)$ and $f_{ij} > 0$ corresponds to the fixed cost of using arc $(i,j)$ to route commodities. For each arc $(i,j) \in A$ and each commodity $k \in K$, $r_{ij}^k >0$
corresponds to the cost of routing one unit of commodity $k$ through arc $(i,j)$.

A MC solution consists of an arc subset $A' \subseteq A$ and, for each commodity $k \in K$, in a flow of value $q^k$ from its origin $o^k$ to its destination $d^k$ with the following requirements: all commodities are only routed through arcs of $A'$, and the total amount of flow routed through each arc $(i,j) \in A'$ does not exceed its capacity $c_{ij}$. The solution cost is the sum of the fixed costs over the arcs of $A'$ plus the routing cost, the latter being the sum over all arcs $(i,j) \in A$ and all commodities $k \in K$ of the unitary routing cost $r_{ij}^k$ multiplied by the amount of flow of $k$ routed through $(i,j)$.

\subsection{MILP formulation}


A standard model for the MC problem~\citep{Gendron1999} introduces two sets of variables: the continuous flow variables $x_{ij}^k$ representing the amount of commodity $k$ that is routed through arc $(i,j)$ and the binary design variables $y_{ij}$ representing whether or not arc $(i,j)$ is used to route commodities.
Denoting respectively by $N^+_i=\{j \in N \mid (i,j)\in A\}$ and $N_i^-=\{j \in N \mid (j, i)\in A\}$ the sets of forward and backward neighbors of a vertex $i \in N$, the MC problem can be modeled as follows:
\begin{subequations}\label{eq:MC_MILP}
\begin{align}
	\label{znd:obj}
	&\min_{\bm{x}, \bm{y}}   \sum_{(i,j) \in A} \left(f_{ij}y_{ij} + \sum_{k \in K} r^k_{ij}x_{ij}^k\right)  \span\span                                       \\
	\label{znd:flowCons}  & \sum_{j\in N^+_i}x_{ij}^k - \sum_{j\in N^-_i}x_{ji}^k=b^k_i & & \!\!\forall i \in N, \forall k \in K\\
	\label{znd:capacity}      & \sum_{k \in K}x_{ij}^k\leq c_{ij}y_{ij},            &  &\!\! \forall (i,j)\in A                  \\
    \label{znd:xij0}  &  x_{ij}^k = 0 && \!\!\!\!\!\begin{array}{l}\forall k \in K, \forall (i,j) \in A  \\ \text{ s.t. }i=d^k \text{ or } j=o^k\end{array}\\
	\label{znd:bound}         & 0\leq x^k_{ij}\leq q^k                              &  & \!\!\forall (i,j)\in A, \forall k \in K \\
	\label{znd:binary}        & y_{ij} \in \{0,1\},                              &     & \!\!\forall (i,j)\in A
\end{align}
\end{subequations}
where
\begin{equation*}
	b^k_i=\left\{ \begin{array}{cc}
		q^k  & \mbox{ if } i = o^k, \\
		-q^k & \mbox{ if } i = d^k, \\
		0    & \mbox{ otherwise.}
	\end{array}\right.
\end{equation*}

The objective function \eqref{znd:obj} minimizes the sum of the routing and fixed costs.  Equations \eqref{znd:flowCons} are the flow conservation constraints that properly define the flow of each commodity through the graph.  Constraints \eqref{znd:capacity} are the capacity constraints ensuring that the total amount of flow routed through each arc does not exceed its capacity or is zero if the arc is not used to route commodities. Equations~\eqref{znd:xij0} ensure that a commodity is not routed on an arc entering its origin or leaving its destination. 
Finally inequalities \eqref{znd:bound} are the bounds for the $x$ variables and inequalities \eqref{znd:binary} are the integer constraints for the design variables.


\subsection{Lagrangian Knapsack Relaxation}

A standard way to obtain good bounds for the MC problem is to solve the \LR{} obtained by dualizing the flow conservation constraints \eqref{znd:flowCons} in formulation  \eqref{znd:obj}-\eqref{znd:binary}.
Let $\pi_i^k$ be the \LM{} associated with node $i \in N$ and commodity $k \in K$.
Dualizing the flow conservation constraints gives the following relaxed Lagrangian problem $LR(\bm{\pi})$\footnote{Since the dualized constraints are equations, $\vpi$ have no sign constraints.}:
$$ \everymath={\displaystyle}
	\begin{array}{l} \min_{(\bm{x},\bm{y}) \text{ satisfies }\eqref{znd:capacity}-\eqref{znd:binary}}    \sum_{(i,j) \in A} \left(f_{ij}y_{ij} + \sum_{k \in K} r^k_{ij}x_{ij}^k\right)\\+\sum_{k \in K}\sum_{i \in N}\pi^k_i\left(b^k_i -\sum_{j\in N^+_i}x_{ij}^k + \sum_{j\in N^-_i}x_{ji}^k\right)\end{array}
$$
Rearranging the terms in the objective function and observing that the relaxed Lagrangian problem is decomposed by arcs, we obtain a subproblem for each arc $(i,j)\in A$ of the form:
\begin{subequations}\label{eq:LRij_MC}
\begin{align}
(LR_{ij}(\bm{\pi}))\quad & \min_{\bm{x}, \bm{y}} f_{ij}y_{ij} + \sum_{k \in K_{ij}} w_{ij}^k x_{ij}^k  \span\span                       \\
& \sum_{k \in K_{ij}}x_{ij}^k\leq c_{ij}y_{ij} \label{LR_MC:capacity}      \\
& 0\leq x^k_{ij}\leq q^k & \forall k \in K_{ij} \label{LR:MC:xBounds}\\
& y_{ij} \in \{0,1\}     & \label{LR_MC:yBinary}
\end{align}
\end{subequations}
where $w_{ij}^k = r_{ij}^k-\pi_{i}^k+\pi_{j}^k$ and $K_{ij} = \{k \in K \mid j \neq o^k \mbox{ and } i \neq d^k\}$ is the set of commodities that may be routed through arc $(i,j)$.

For each $(i,j) \in A$, $LR_{ij}(\bm{\pi})$ is a MILP with only one binary variable. If $y_{ij} = 0$, then, by \eqref{LR_MC:capacity} and \eqref{LR:MC:xBounds}, $x_{ij}^k = 0$ for all $k \in K_{ij}$. If $y_{ij} = 1$, the problem reduces to a continuous knapsack problem. An optimal solution is obtained by ordering the commodities of $K_{ij}$ with respect to decreasing values $w_{ij}^k$ and setting for each variable $x_{ij}^k$ the value $\max\{\min\{q^k, c_{ij} - \sum_{k \in K(k)} q^k\}, 0\}$ where $K(k)$ denotes the set of commodities that preceded $k$ in the order. This step can be done in $O(|K_{ij}|)$ if one computes $x_{ij}^k$ following the computed order. Hence, the complexity of the continuous knapsack problem is $O(|K_{ij}|\log(|K_{ij}|))$. The solution of $LR_{ij}(\bm{\pi})$ is the minimum between the cost of the continuous knapsack problem and $\mathbf{0}$.

Lagrangian duality implies that 
\[LR(\bm{\pi})=\sum_{(i,j) \in A} LR_{ij}(\bm{\pi})+ \sum_{i \in N}\sum_{k \in K} \pi_i^k b^k_i
\] 
is a lower bound for the MC problem and the best one is obtained by solving the following Lagrangian dual problem:
\begin{equation*}
	(LD) \quad  \max_{\bm{\pi} \in \mathbb{R}^{N \times K}} LR(\bm{\pi})
\end{equation*}

\section{Generalized Assignment Problem}
\label{app:GA}

A GA instance is defined by a set $I$ of items and a set $J$ of bins.
Each bin $j$ is associated with a certain capacity $c_j$.
For each item $i \in I$ and each bin $j \in J$, $p_{ij}$ is the profit of assigning item $i$ to bin $j$, and $w_{ij}$ is the weight of item $i$ inside bin $j$.

Considering a binary variable $x_{ij}$ for each item and each bin that is equal to one if and only if item $i$ is assigned to bin $j$, the GA problem can be formulated as:
\begin{subequations}\label{eq:GA_MILP}
\begin{align}
& \max_{\bm{x}}  \sum_{i \in I} \sum_{j \in J} p_{ij}x_{ij} \label{obj:GA}\\
& \sum_{j \in J} x_{ij}  \leq 1 && \forall i \in I \label{sa_constr:GA} \\
& \sum_{i \in I} w_{ij} x_{ij}  \leq c_j && \forall j \in J \label{cap_constr:GA}\\
& x_{ij} \in \{0,1\} && \forall i \in I ,\; \forall j\in J .\label{int_constr:GA}
\end{align}
\end{subequations}

The objective function \eqref{obj:GA} maximizes the total profit.
Inequalities \eqref{sa_constr:GA} assert that each item is contained in no more than one bin.
Inequalities \eqref{cap_constr:GA} ensure that the sum of the weights of the items assigned to a bin does not exceed its capacity.
Finally, constraints \eqref{int_constr:GA} assure the integrality of the variables.

\subsection{Lagrangian Relaxation}

A \LR{} of the GA problem is obtained by dualizing \eqref{sa_constr:GA}.
For $i \in I$, let $\pi_{i} \geq 0$ be the \LM{} of inequality \eqref{sa_constr:GA} associated with item $i$.
For each bin $j$ the subproblem becomes:
\begin{align*}
(LR_j(\bm{\pi})) \quad & \max_{\bm{x}} \sum_{i \in I}\sum_{j \in J} (p_{ij}-\pi_{i}) x_{ij} \span \span\\
& \sum_{i \in I}w_{ij}x_{ij} \leq c_j\\
& x_{ij} \in \{0,1\} && \forall i \in I
\end{align*}
It corresponds to an integer knapsack with $|I|$ binary variables. For $\bm{\pi}\ge \bm{0}$, the Lagrangian bound $LR(\bm{\pi})$ is:
\begin{equation*}
    LR(\bm{\pi}) = \sum_{j \in J}LR_j(\bm{\pi})+\sum_{i \in I}\pi_i.
\end{equation*}

The Lagrangian dual can then be written as:
\begin{equation*}
    \min_{\bm{\pi} \in \sR^{|I|}_{\geq 0}}LR(\bm{\pi})
\end{equation*}

\section{Hyperparameters and Implementation Details}\label{sec:hyper_parameters_and_implementation}

In this section, we provide technical details about the implementation.
We also specify the hyperparameters of the NN architecture (regarding layers, layer sizes, and activation functions) for the architecture, define the features (hand-crafted), specify the hyperparameters for the optimizer used for the training, and the GPU and CPU specifics of the machines used for the numerical experiments.

\paragraph{Loss functions}

To weigh the contributions of the different iterations, we use $\gamma=0.999$.
In other amortized optimization works, focusing on learning optimizers, as~\citep{NIPS2016_fb875828}, this parameter is taken equal to $1$.

\begin{table}[h]
\centering
\caption{Best average initialization values for the $\eta_0$ parameter over the test set for all methods that require a grid search. The values $\{10^{4},\, 10^{3},\, 10^{2},\, 10^{1},\, 10^{0},\, 10^{-1}\}$ are evaluated, and for each fixed maximum number of iterations ($10$, $25$, $50$, and $100$), the value yielding the lowest average GAP is selected.}
    \label{tab:grid_search}
\begin{tabular}{@{}ll rrrr@{}}
{Dataset} & {Methods} & 
10 iter. & 
25 iter. & 
50 iter. & 
100 iter. \\
\midrule
\multirow{6}{*}{{\tiny\textsc{MC-Sml-40}}} 
& Bundle h. & 10 & 10 & 10 & 10 \\
& Bundle b. & 100 & 10 & 10 & 10 \\
& Bundle s. & 100 & 10 & 10 & 100 \\
& Bundle c. & 100 & 100 & 100 & 100 \\
& Descent & 10 & 10 & 10 & 10 \\
& Adam & 1 & 1 & 1 & 1 \\
\midrule
\multirow{6}{*}{{\tiny\textsc{MC-Sml-Var}}} 
& Bundle h. & 10 & 100 & 100 & 10 \\
& Bundle b. & 100 & 100 & 100 & 100 \\
& Bundle s. & 100 & 100 & 100 & 100 \\
& Bundle c. & 100 & 100 & 100 & 100 \\
& Descent & 100 & 10 & 10 & 10 \\
& Adam & 1 & 1 & 1 & 1 \\
\midrule
\multirow{6}{*}{{\tiny\textsc{MC-Big-40}}} 
& Bundle h. & 10 & 10 & 10 & 1 \\
& Bundle b. & 10 & 10 & 10 & 10 \\
& Bundle s. & 10 & 10 & 10 & 10 \\
& Bundle c. & 10 & 10 & 10 & 10 \\
& Descent & 10 & 10 & 10 & 10\\
& Adam & 1 & 1 & 1 & 1 \\
\midrule
\multirow{6}{*}{{\tiny\textsc{MC-Big-Var}}}
& Bundle h. & 100 & 10 & 10 & 10 \\
& Bundle b. & 100 & 100 & 10 & 10 \\
& Bundle s. & 100 & 100 & 10 & 10 \\
& Bundle c. & 100 & 100 & 100 & 100 \\
& Descent & 10 & 10 & 10 & 10\\
& Adam & 1 & 1 & 1 & 1 \\
\midrule
\multirow{6}{*}{{\tiny\textsc{GA-10-100}}} 
& Bundle h. & 1000 & 1000 & 1000 & 1000 \\
& Bundle b. & 1000 & 1000 & 1000 & 1000 \\
& Bundle s. & 1000 & 1000 & 1000 & 1000 \\
& Bundle c. & 1000 & 1000 & 1000 & 1000 \\
& Descent & 100 & 100 & 100 & 100 \\
& Adam & 100 & 10 & 10 & 10 \\
\midrule
\multirow{6}{*}{{\tiny\textsc{GA-20-400}}} 
& Bundle h. & 1000 & 1000 & 1000 & 1000 \\
& Bundle b. & 1000 & 1000 & 1000 & 1000 \\
& Bundle s. & 1000 & 1000 & 1000 & 1000 \\
& Bundle c. & 1000 & 1000 & 1000 & 1000 \\
& Descent & 1000 & 100 & 100 & 100 \\
& Adam & 100 & 10 & 10 & 10 \\
\bottomrule
\end{tabular}
\end{table}

\paragraph{Model Architecture - Bundle Network}

For all datasets, we consider the same architecture.
A first LSTM that goes from the feature space to predict $6$ vectors of size $128$.
These are the means and variances used to sample three vectors that are the hidden representation of the $\eta$ parameter, the one for the key for the last found component, and the one for the query for the current iteration. These vectors are given as inputs to three different decoders, each one composed of a hidden layer of size $8 \cdot 128$. The output of these decoders is respectively a scalar $\eta_t$, used as step size, and two vectors of size $128$ to represent the key and the query of the current iteration.
The main framework is presented with a sample mechanism, but we found that it was detrimental when using softmax for function $\psi$, and we refer to Section~\ref{sec:sampling} for further details. Still, sampling improves performance when used in conjuction with sparsemax for $\psi$.

\paragraph{Features}

At iteration $t$, we perform a human-designed feature extraction, chosen to provide similar information to that used by the heuristics $\eta-$strategies. We add further features to represent the component added to the bundle in the associated iteration, and we others to characterize this entry and to compare it with previously added components. Features related to iteration $t$, summarizing the optimization dynamics at the current iteration:
\begin{itemize}
    \item the last step-size $\eta_{t-1}$,
    \item The square norm of the search direction $||\vw^{(t-1)}||^2_2$,
    \item The square norm of the search direction weighted with $\eta_t$, that is $\eta_{t-1}||\vw^{(t-1)}||_2^2$. This corresponds to the quadratic part in the DMP objective function,
    \item The linear part in the DMP: $\sum_{j=1}^{t-1}\alpha_j\theta_{j}^{(t-1)}$,
    \item a boolean value to see if the quadratic part is bigger than the linear part $||\vw^{(t-1)}||^2_2>\sum_{j=1}^{t-1}\alpha_j\theta_{j}^{(t-1)}$,
    \item a boolean to see if the quadratic part, rescaled by a \textit{big} $\eta^*=10000$, is greater than the linear part: $\eta^*||\vw^{(t-1)}||^2_2>\sum_{j=1}^{t-1}\alpha_j\theta_{j}^{(t-1)}$,
    \item The iteration counter $t$.
\end{itemize}
Features describing the last trial and stabilization points. These features capture information about the most recent trial point $\vpi_t$ and stabilization point $\bar{\vpi}_t$:
\begin{itemize}
    \item The objective value in the last trial point $\phi(\vpi_{t})$ and in the last stabilization point $\phi(\bar{\vpi}_{t})$,
    \item The linearization error in the last stabilization point $\bar{\alpha}_{t}$ and in the last trial point $\alpha_{t}$,
    \item The square norm of the last trial point $||\vpi_{t}||_2$ of the last stabilization point $||\bar{\vpi}_{t}||_2$ and of the gradient in the stabilization point $||\bar{\vg}_{t}||_2$,
    \item the square norm, the mean, the variance, the minimum, and the maximum of the vector $\vg_{t}$,
    \item the square norm, the mean, the variance, the minimum, and the maximum of the vector  $\vpi_{t}$.
\end{itemize}
Features comparing the last inserted component with the ones already contained in the bundle:
\begin{itemize}
    \item The minimum and the maximum of the scalar product of the gradients in the bundle $\min_j(\vg_{t}^{\top}\vg_{j})$ and $\max_j(\vg_{t}^{\top}\vg_{j})$,
    \item The minimum and the maximum of the scalar product of the points in the bundle $\min_j(\vpi_{t}^{\top}\vpi_{j})$ and $\max_j(\vpi_{t}^{\top}\vpi_{j})$,
    \item The scalar product of the last inserted gradient and the last search direction $\vg_{t}^{\top} \vw^{(t-1)}$.
\end{itemize}

\paragraph{Optimiser Specifications}

We use Adam as optimizer, with a learning rate $0.00001$, a Clip Norm (to 5), and exponential decay $0.9$.

\paragraph{GPU specifics}

For the training on the datasets, we use Quadro RTX 5000 accelerators with 8Gb of RAM.
We test performance on the same machine.
All the variants, except for Bundle Network, are CPU-only based. The experiments are done on the same machine with QEMU Virtual CPU version 2.5+.

\section{Dataset generalization sparsemax and sampling}

In this appendix, we provide some further details on the model that uses sparsemax and the sampling mechanism, related to cross-dataset generalization properties, similarly to what was done in Section \ref{sec:cross_testing} for the softmax without sampling.

\begin{table}[h!]
\caption{Gap of the model tested on MC-Sml-40 trained on different datasets. The associated training dataset is displayed as first element of each row of the table. The last row shows the performances of the best baseline obtained with grid-search.}
\begin{tabular}{@{}l c rr c rr c rr c rr@{}}
\toprule
\multirow{2}{*}{Training Dataset} & & \multicolumn{2}{c}{10 iter.} & 
& \multicolumn{2}{c}{25 iter.} & 
& \multicolumn{2}{c}{50 iter.} & 
& \multicolumn{2}{c}{100 iter.} \\
\cmidrule{3-4} \cmidrule{6-7} \cmidrule{9-10} \cmidrule{12-13}
& & GAP & time & 
& GAP & time & 
& GAP & time & 
& GAP & time \\
\midrule
\textbf{MC-Sml-40} & & 10.51 & 0.126 & & 2.07 & 0.214 & & 0.60 & 0.363 & & \textbf{0.18} & 0.662\\
MC-Sml-Var & & 9.87 & 0.128 & & \textbf{1.75} & 0.218 & & \textbf{0.54} & 0.371 & & \textbf{0.18} & 0.683\\
MC-Big-40 & & 21.36 & 0.140 & & 5.23 & 0.239 & & 1.78 & 0.41 & & 0.52 & 0.755\\
MC-Big-Var & & \textbf{9.82} & 0.128 & & 1.99 & 0.221 & & 0.60 & 0.374 & & 0.19 & 0.683\\
\midrule
GA-10-100 & & 69.18 & 0.128 & & 69.18 & 0.216 & & 69.18 & 0.366 & & 69.18 & 0.667\\
GA-20-400 & & 96.08 & 0.128 & & 96.08 & 0.218 & & 96.08 & 0.368 & & 96.08 & 0.673\\
\midrule
Best Baseline & & 11.96 & 0.048 & & 4.16 & 0.128 & & 1.30 & 0.293 & & 0.31 & 0.816\\
\bottomrule
     \end{tabular}
\end{table}

\begin{table}[h!]
\caption{Gap of the model tested on MC-Sml-Var and trained on different datasets. The associated training dataset is displayed as first element of each row of the table. The last row shows the performances of the best baseline obtained with grid-search.}
\begin{tabular}{@{}l c rr c rr c rr c rr@{}}
\toprule
\multirow{2}{*}{Training Dataset} & & \multicolumn{2}{c}{10 iter.} & 
& \multicolumn{2}{c}{25 iter.} & 
& \multicolumn{2}{c}{50 iter.} & 
& \multicolumn{2}{c}{100 iter.} \\
\cmidrule{3-4} \cmidrule{6-7} \cmidrule{9-10} \cmidrule{12-13}
& & GAP & time & 
& GAP & time & 
& GAP & time & 
& GAP & time \\
\midrule
MC-Sml-40 & & 22.9 & 0.161 & & 6.62 & 0.301 & & 2.6 & 0.54 & & 1.02 & 1.027\\
\textbf{MC-Sml-Var} & & {14.02} & 0.158 & & 3.80 & 0.298 & & 1.56 & 0.531 & & 0.64 & 1.007\\
MC-Big-40 & & 39.34 & 0.171 & & 15.12 & 0.327 & & 6.7 & 0.59 & & 2.97 & 1.119\\
MC-Big-Var & & 14.27 & 0.160 & & \textbf{3.64} & 0.302 & & \textbf{1.27} & 0.545 & & \textbf{0.49} & 1.042\\
\midrule
GA-10-100 & & 59.75 & 0.160 & & 59.75 & 0.299 & & 59.75 & 0.540 & & 59.75 & 1.031\\
GA-20-400 & & 66.35 & 0.163 & & 66.35 & 0.307 & & 66.35 & 0.552 & & 66.35 & 1.052\\
\midrule
Best Baseline & & \textbf{12.79} & 0.108 &  & 4.33 & 0.29 & & 1.66  &  0.635 & & 0.58 & 1.503\\
\bottomrule
\end{tabular}
\end{table}

\begin{table}[h!]
\caption{Gap of the model tested on MC-Big-40 and trained on different datasets. The associated training dataset is displayed as first element of each row of the table. The last row shows the performances of the best baseline obtained with grid-search.}
\begin{tabular}{@{}l c rr c rr c rr c rr@{}}
\toprule
\multirow{2}{*}{Training Dataset} & & \multicolumn{2}{c}{10 iter.} & 
& \multicolumn{2}{c}{25 iter.} & 
& \multicolumn{2}{c}{50 iter.} & 
& \multicolumn{2}{c}{100 iter.} \\
\cmidrule{3-4} \cmidrule{6-7} \cmidrule{9-10} \cmidrule{12-13}
& & GAP & time & 
& GAP & time & 
& GAP & time & 
& GAP & time \\
\midrule
MC-Sml-40 & & 41.77 & 0.154 & & 7.62 & 0.276 & & 1.89 & 0.484 & & 0.67 & 0.915\\
MC-Sml-Var & & 91.76 & 0.152 & & 63.78 & 0.279 & & 37.85 & 0.487 & & 8.22 & 0.919\\
\textbf{MC-Big-40} & & \textbf{13.27} & 0.150 & & \textbf{3.50} & 0.272 & & \textbf{1.40} & 0.479 & & \textbf{0.51} & 0.905\\
MC-Big-Var & & 44.89 & 0.154 & & 7.19 & 0.280 & & 2.17 & 0.492 & & 1.09 & 0.931\\
\midrule
GA-10-100 & & 99.8 & 0.140 & & 99.8 & 0.252 & & 99.8 & 0.445 & & 99.8 & 0.84\\
GA-20-400 & & 100.0 & 0.156 & & 100.0 & 0.281 & & 100.0 & 0.493 & & 100.0 & 0.931\\
\midrule
Best Baseline & & 19.34 & 0.076 & & 5.64 & 0.203 & & 2.08 &  0.447 & & 0.78 & 1.110\\
\bottomrule
     \end{tabular}
\end{table}

\begin{table}[h!]
\caption{Gap of the model tested on MC-Big-Var and trained on different datasets. The associated training dataset is displayed as first element of each row of the table. The last row shows the performances of the best baseline obtained with grid-search.}
\begin{tabular}{@{}l c rr c rr c rr c rr@{}}
\toprule
\multirow{2}{*}{Training Dataset} & & \multicolumn{2}{c}{10 iter.} & 
& \multicolumn{2}{c}{25 iter.} & 
& \multicolumn{2}{c}{50 iter.} & 
& \multicolumn{2}{c}{100 iter.} \\
\cmidrule{3-4} \cmidrule{6-7} \cmidrule{9-10} \cmidrule{12-13}
& & GAP & time & 
& GAP & time & 
& GAP & time & 
& GAP & time \\
\midrule
MC-Sml-40 & & 24.9 & 0.162 & & 6.58 & 0.307 & & 2.35 & 0.555 & & 0.89 & 1.057\\
MC-Sml-Var & & 29.88 & 0.252 & & 14.63 & 0.46 & & 7.67 & 0.822 & & 1.85 & 1.546\\
MC-Big-40 & & 32.67 & 0.174 & & 12.04 & 0.33 & & 5.31 & 0.598 & & 2.30 & 1.141\\
\textbf{MC-Big-Var} & & {18.94} & 0.162 & & \textbf{4.24} & 0.305 & & \textbf{1.47} & 0.548 & & \textbf{0.61} & 1.045\\
\midrule
GA-10-100 & & 67.42 & 0.165 & & 67.42 & 0.31 & & 67.42 & 0.561 & & 67.42 & 1.069\\
GA-20-400 & & 74.69 & 0.163 & & 74.69 & 0.306 & & 74.69 & 0.551 & & 74.69 & 1.051\\
\midrule
Best Baseline & & \textbf{17.98} & 0.112 & & 6.34 & 0.289 &  & 2.66 & 0.627 & & 0.93 & 1.477\\
\bottomrule     \end{tabular}
\end{table}

\begin{table}[h!]
\caption{Gap of the model tested on GA-10-100 and trained on different datasets. The associated training dataset is displayed as first element of each row of the table. The last row shows the performances of the best baseline obtained with grid-search.}
\begin{tabular}{@{}l c rr c rr c rr c rr@{}}
\toprule
\multirow{2}{*}{Training Dataset} & & \multicolumn{2}{c}{10 iter.} & 
& \multicolumn{2}{c}{25 iter.} & 
& \multicolumn{2}{c}{50 iter.} & 
& \multicolumn{2}{c}{100 iter.} \\
\cmidrule{3-4} \cmidrule{6-7} \cmidrule{9-10} \cmidrule{12-13}
& & GAP & time & 
& GAP & time & 
& GAP & time & 
& GAP & time \\
\midrule
MC-Sml-40 & & 2.55 & 0.123 & & 1.79 & 0.208 & & 0.85 & 0.35 & & 0.33 & 0.636\\
MC-Sml-Var & & 2.42 & 0.121 & & 1.4 & 0.202 & & 0.65 & 0.341 & & 0.56 & 0.613\\
MC-Big-40 & & 2.67 & 0.119 & & 2.12 & 0.197 & & 1.4 & 0.334 & & 0.69 & 0.603\\
MC-Big-Var & & 2.5 & 0.120 & & 1.61 & 0.203 & & 0.9 & 0.339 & & 0.59 & 0.609\\
\midrule
\textbf{GA-10-100} & & 0.92 & 0.122 & & 0.70 & 0.201 & & 0.68 & 0.34 & & 0.68 & 0.609\\
GA-20-400 & & \textbf{0.37} & 0.122 & & 0.11 & 0.207 & & 0.06 & 0.352 & & 0.05 & 0.636\\
\midrule
Best Baseline & & {0.1893} & 0.071 & & \textbf{0.0156} & 0.211 & & \textbf{0.0014} & 0.601 & & \textbf{0.0006} & 1.596\\
\bottomrule
     \end{tabular}
\end{table}

\begin{table}[h!]
\caption{Gap of the model tested on GA-20-400 and trained on different datasets. The associated training dataset is displayed as first element of each row of the table. The last row shows the performances of the best baseline obtained with grid-search.}
\begin{tabular}{@{}l c rr c rr c rr c rr@{}}
\toprule
\multirow{2}{*}{Training Dataset} & & \multicolumn{2}{c}{10 iter.} & 
& \multicolumn{2}{c}{25 iter.} & 
& \multicolumn{2}{c}{50 iter.} & 
& \multicolumn{2}{c}{100 iter.} \\
\cmidrule{3-4} \cmidrule{6-7} \cmidrule{9-10} \cmidrule{12-13}
& & GAP & time & 
& GAP & time & 
& GAP & time & 
& GAP & time \\
\midrule
MC-Sml-40 & & 2.35 & 0.253 & & 1.9 & 0.536 & & 1.14 & 1.004 & & 0.39 & 1.943\\
MC-Sml-Var & & 2.27 & 0.238 & & 1.67 & 0.504 & & 0.76 & 0.937 & & 0.5 & 1.805\\
MC-Big-40 & & 2.42 & 0.237 & & 2.08 & 0.504 & & 1.58 & 0.949 & & 0.88 & 1.836\\
MC-Big-Var & & 2.32 & 0.237 & & 1.79 & 0.503 & & 1.01 & 0.949 & & 0.6 & 1.84\\
\midrule
GA-10-100 & & 1.17 & 0.252 & & 1.16 & 0.537 & & 1.16 & 1.005 & & 1.16 & 1.943\\
\textbf{GA-20-400} & & {0.4} & 0.218 & & {0.13} & 0.461 & & {0.07} & 0.86 & & {0.06} & 1.651\\
\midrule
Best Baseline & & \textbf{0.1129} & 0.433 & & \textbf{0.0193} & 1.459 & & \textbf{0.0048} & 4.256 & & \textbf{0.0020} & 14.459\\
\bottomrule
     \end{tabular}
\end{table}

\end{appendices}


\end{document}